\documentclass[11pt]{article}
\usepackage{amssymb}
\usepackage{amsmath}
\usepackage{verbatim}
\usepackage{enumerate}
\usepackage{hyperref}

\usepackage{amscd,enumitem}
\usepackage{amsthm,mathtools}
\usepackage{indentfirst}
\hypersetup{
  colorlinks   = true, 
  urlcolor     = blue, 
  linkcolor    = blue, 
  citecolor   = red 
}

\newtheorem{theorem}{Theorem}[section]
\newtheorem{lemma}[theorem]{Lemma}
\newtheorem{corollary}[theorem]{Corollary}
\newtheorem{proposition}[theorem]{Proposition}

\newtheorem{remark}[theorem]{Remark}

\theoremstyle{definition}

\DeclareMathOperator{\id}{id}
\DeclareMathOperator{\vac}{\mathbf{1}}

\DeclareMathOperator{\SL}{SL}

\DeclareMathOperator{\tr}{tr}
\DeclareMathOperator{\str}{str}
\DeclareMathOperator{\Aut}{Aut}
\DeclareMathOperator{\End}{End}
\DeclareMathOperator{\wt}{wt}
\DeclareMathOperator{\PP}{\mathbb{P}}
\DeclareMathOperator{\slthat}{\widehat{\frak{sl}}_2}
\DeclareMathOperator{\slt}{\frak{sl}_2}

\numberwithin{equation}{section}

\newcommand{\Qed}{\quad \square}

\newcommand{\tpi}{2\pi i}

\newcommand{\trq}[2]{\tr_{#1}{#2}q^{L(0)}}

\newcommand{\half}{\frac{1}{2}}

\newcommand{\Z}{\mathbb{Z}}
\newcommand{\C}{\mathbb{C}}

\newcommand{\N}{\mathbb{N}}

\newcommand{\HH}{\mathbb{H}}

\let\OLDthebibliography\thebibliography
\renewcommand\thebibliography[1]{
  \OLDthebibliography{#1}
  \setlength{\parskip}{0pt}
  \setlength{\itemsep}{0pt plus 0.3ex}
}

\begin{document}

\title{Zhu reduction for Jacobi $n$-point functions and applications}
\author{
Kathrin Bringmann\thanks{Institute of Mathematics, University of Cologne. E-mail: kbringma@math.uni-koeln.de. The research is supported by the Alfried Krupp Prize for Young University
Teachers  of  the  Krupp  foundation  and  the  research  leading  to  these  results  receives  funding
from the European Research Council under the European Union's Seventh Framework Programme
(FP/2007-2013) / ERC Grant agreement n.  335220 - AQSER. },\quad
Matthew Krauel\thanks{Department of Mathematics and Statistics, California State University, Sacramento. E-mail: krauel@csus.edu.  Research for this project was supported by the European Research Council (ERC) Grant agreement n. 335220 - AQSER, and also by a research visit to the Max Planck Institute for Mathematics in Bonn.},\quad
Michael P.\ Tuite\thanks{School of Mathematics, Statistics and Applied Mathematics, NUI Galway. E-mail: michael.tuite@nuigalway.ie.}
}

\date{}
\maketitle

\begin{abstract}
\noindent We establish precise Zhu reduction formulas for Jacobi $n$-point functions which show the absence of any possible poles arising in these formulas. We then exploit this to produce results concerning the structure of strongly regular vertex operator algebras, and also to motivate new differential operators acting on Jacobi forms. Finally, we apply the reduction formulas to the Fermion model in order to create polynomials of quasi-Jacobi forms which are Jacobi forms.
\end{abstract}

\date{}
\maketitle


\section{Introduction}\label{Sec:Intro}

One of the most foundational works in the theory of Vertex Operator Algebras (VOAs) (e.g. \cite{FHL,  Kac-VOA, LL}) is Zhu's study of elliptic $n$-point functions \cite{Zhu}.
Zhu developed important reduction (or recursion) formulas which allow, among other things, $n$-point functions to be written as linear combinations of $(n-1)$-point functions with coefficients that are quasi-modular forms. The fact that quasi-modular forms are holomorphic on the complex upper half-plane, $\mathbb{H}$, then helps to show that all $n$-point functions are holomorphic on this domain in many cases. Meanwhile, a desire to study Jacobi $n$-point functions, involving an additional variable $z\in\C$, leads us to consider similar reduction formulas for these generalizations. The coefficients that arise in these reductions, however, are in some cases so-called quasi-Jacobi forms with simple poles in the variable $z\in\C$. For this reason, it is conceivable that $1$-point functions that are descendants of the vacuum element may have poles. One aim of this paper is to explain how the Zhu reduction formulas for Jacobi $n$-point functions do not in fact introduce such poles.
This follows from an analysis of Zhu reduction in the neighborhood of all possible poles i.e., not just near $z=0$.

Zhu's work has been extremely useful in progressing the study of single variable $n$-point functions (for example, \cite{DLM-orbifold, MTZ, Miy-C2}). Meanwhile, although Jacobi $n$-point functions often serve advantages (for example, unlike their single-valued brethren they discern the difference between inequivalent irreducible modules for VOAs associated with affine Lie algebras), there has yet to be a complete analysis of the Zhu reduction formulas. In particular, existing formulas either avoid such poles \cite{KrauelMasonII} or apply to an example where an alternative approach to explaining the lack of poles may be taken \cite{GK-differential}. After introducing the relevant functions in Section \ref{Sec:Quasi}, we turn to establishing the Zhu reduction formulas for Jacobi $n$-point functions in Section \ref{Sec:Zhu}. We also include here an alternative approach using the shifted theories for VOAs (see \cite{DM-shifted, Li} for discussions on such theories).
As a corollary to Propositions \ref{prop:Zhured}--\ref{prop:apnpt0} below, we obtain the following theorem.
\begin{theorem}\label{thm:intro1}
 A Jacobi $n$-point function for a VOA $V$ does not contain poles in $\mathbb{C}\times \mathbb{H}$ if the $(n-1)$-point functions do not contain poles for any $n-1$ many vectors in $V$.
\end{theorem}

The reasons why poles exist in the reduction formula coefficients and yet not in the associated $n$-point functions are interesting in their own right, and quite exploitable. Indeed, the bulk of this paper examines this process in more detail. At its core, the possible poles must either never exist (i.e., there are no elements in the VOA which produce the quasi-Jacobi forms giving rise to the poles), or the poles must correspond to zeros of the partition function. Gaberdiel and Keller \cite{GK-differential} used the Zhu reduction formulas for Jacobi $n$-point functions (or elliptic genus) and the fact that no poles arise in the $N=2$ superconformal field theories to create new differential operators of Jacobi forms of different (higher) degrees. They also highlighted the use of this for investigating extremal $N=2$ superconformal field theories.

In Section \ref{Sec:DiffOpsOnJacobi} we study a family of differential operators $\mathcal{M}_{k,\alpha}$ defined for $k\in \mathbb{N}$ and $\alpha \in \mathbb{Z}\setminus \{0\}$. For certain $k$ and $\alpha$ these operators collapse to those studied in \cite{GK-differential} and for other values to those considered in \cite{Oberdieck-Serre}. However, some subtle additional cases are included here. Along with showing certain coefficients of functions under the image of this form are nonzero (see Lemma \ref{lem:JacobiCoeffs2}), we also establish Lemma \ref{lem:LemPreserve} which, in a simplified version, can be paraphrased as follows.
\begin{lemma}
Suppose $\alpha \in \mathbb{Z}\setminus \{0\}$, $k,m \in \mathbb{N}_0$, and $\phi$ is a weak Jacobi form of weight $k$ and index $m$ (with a possible multiplier system). Then $\mathcal{M}_{k,\alpha} (\phi)$ transforms like a Jacobi form of weight $k+2$ and index $m$ (with the same multiplier system). Additionally, if $\alpha = \pm 1, \pm 2$, then $\mathcal{M}_{k,\alpha} (\phi)$ is holomorphic for either $k$ even or odd if certain conditions on the multiplier system are satisfied.
\end{lemma}

The operators studied in Section \ref{Sec:DiffOpsOnJacobi} are motivated by applying the Zhu reduction formulas to strongly regular VOAs. This analysis is performed in Section \ref{Sec:DiffOpsOnVOAs}. We highlight the fact that an additional Lie algebra structure contained in the strongly regular VOAs is what gives rise to the new differential operators. This is described further in Section \ref{Sec:DiffOpsOnVOAs}. One could also consider higher degree differential operators here, much as in the same way as Gaberdiel and Keller do in \cite{GK-differential} but where one does not have the Lie algebra structure. However, while interesting, this is not pursued in the present paper. Instead, we develop some applications of the existence of the degree $2$ operator and also consider two examples of strongly regular VOAs.

While the Fermion model is not technically a VOA (but rather a vertex operator super algebra (VOSA)), we explain in Section \ref{Sec:Fermion} how we are also able to analyze this VOSA with the Zhu reduction formulas. Among finding degree two differential operators preserving Jacobi forms here that we cannot find for strongly regular VOAs, we also are able to find a degree one operator. We then use the developed theory to find and study polynomials of quasi-Jacobi forms which are Jacobi forms.

Finally, we mention a few (there are many other) instances in the literature where our work intersects. Quasi-Jacobi forms play a significant role in vertex algebra theory in the work \cite{HE-characters}, where the characters of topological $N=2$ vertex algebras are studied and found to be Jacobi forms. Additionally, the study of Gromov-Witten potentials \cite{Kaw} provides another vantage point of quasi-Jacobi forms in a related field. Calculating elliptic genera, which are closely related to the Jacobi partition (or $0$-point) functions, considered here, for Landau-Ginzburg orbifolds can be found in \cite{KYY}. Additional work dealing with elliptic genera \cite{Libgober-Elliptic} explores quasi-Jacobi forms in more depth, and is used often here.

%


\section{Jacobi and quasi-Jacobi forms}\label{Sec:Quasi}

\subsection{Basic definitions}

We start by recalling classical Jacobi forms, the reader is referred to \cite{EZ} for good background material. Let $k,m\in\N_0$. A {\it holomorphic Jacobi form} of weight $k$ and index $m$ on $\SL_2(\Z)$ with rational multiplier $\chi$ (a rational character for a one dimensional representation of the Jacobi group $\SL(2,\Z)\ltimes \Z^2$  is a holomorphic function $\phi: \C \times \HH \to\C$ which satisfies the following conditions:
\begin{enumerate}[label=\normalfont(\roman*)]
\item We have for $\begin{psmallmatrix}a&b\\c&d\end{psmallmatrix}\in\SL_2(\Z)$ and $(\lambda,\mu)\in\Z^2$
\begin{equation}\label{eq:Jactr}
\phi\Big|_{k,m}\left(\begin{pmatrix}a&b\\c&d\end{pmatrix},(\lambda,\mu)\right)=\chi \left(\begin{pmatrix}a&b\\c&d\end{pmatrix},(\lambda,\mu)\right)\phi
\end{equation}
where for a function $\phi: \C \times \HH \to\C$
\begin{multline*}
\phi\Big|_{k,m}\left(\begin{pmatrix}a&b\\c&d\end{pmatrix},(\lambda,\mu)\right)(z,\tau)\\
:= (c\tau+d)^{-k}e\left(-\frac{cm(z+\lambda\tau+\mu)^2}{c\tau+d}+m\left(\lambda^2\tau+2\lambda z\right)\right)\phi\left(\frac{z+\lambda\tau+\mu}{c\tau+d}, \frac{a\tau+b}{c\tau+d}\right).
\end{multline*}
Here $e(w):=e^{2\pi iw}$.

\item For a multiplier $\chi$, we abbreviate  for $\begin{psmallmatrix}a&b\\c&d\end{psmallmatrix}\in\SL_2(\Z)$ and $(\lambda,\mu)\in\Z^2$
\[
\chi\begin{pmatrix}a&b\\c&d\end{pmatrix} :=\chi\left(\begin{pmatrix}a&b\\c&d\end{pmatrix}, (0, 0)\right),\quad
\chi(\lambda, \mu):=\chi\left(\begin{pmatrix}1&0\\0&1\end{pmatrix},
(\lambda, \mu)\right),
\]
and let $N_1, N_2\in\N$ be uniquely defined by
\[
\chi\left(\begin{matrix} 1&1\\0&1\end{matrix}\right):=e^{2\pi i\frac{a_1}{N_1}}, \quad\chi(0, 1):=e^{2\pi i\frac{a_2}{N_2}},
\]
where $a_j \in\N$ satisfy $\gcd (a_j, N_j)=1$.

In terms of $q:=e(\tau)$ and $\zeta :=e(z)$, the function $\phi$ has a Fourier expansion of the form
 \begin{align}
\phi \left(z ,\tau \right)= \sum_{n\in \N_0+\rho_{1}}
 \sum_{r\in\Z+\rho_{2}\atop{r^2\leq 4nm}}
c(n,r)q^{n}\zeta^{r},
\label{eq:phiFourier}
 \end{align}
where $\rho_{j}:=\frac{a_j}{N_j}\pmod \Z$ with $0\le \rho_j<1$.
\end{enumerate}

If additionally in (ii), $\phi$ satisfies the condition $c(n,r)=0$ if $4mn=r^2$, then $\phi$ is called a {\it Jacobi cusp form}.  If the condition $4mn \geq r^2$ is replaced with the weaker condition $r\in \mathbb{Z}$, then $\phi$ is referred to as a {\it weak Jacobi form}.

Note that the holomorphic Jacobi forms (respectively, Jacobi cusp forms, weak Jacobi forms) of weight $k$ and index $m$ naturally form a $\C$-vector space which we denote by $J_{k,m,\chi}$ (respectively, $J_{k,m, \chi}^{\rm cusp}$, $\widetilde{J}_{k,m, \chi}$). We also consider {\it meromorphic Jacobi forms} which allow poles in the elliptic $z$-variable.

We next consider quasi-Jacobi forms as introduced by Libgober \cite{Libgober-Elliptic}. An {\it almost meromorphic Jacobi form} of weight $k$, index $0$, and depth $(s,t)$ is a meromorphic function in $\C\{q,\zeta\}[z^{-1},\frac{z_2}{\tau_2},\frac{1}{\tau_2}]$ (where $z=z_1+iz_2,\ \tau=\tau_1+i\tau_2$) which satisfies \eqref{eq:Jactr} and has degree at most $s$ and $t$ in $\frac{z_2}{\tau_2}$ and $\frac{1}{\tau_2}$, respectively. A {\it quasi-Jacobi form} of weight $k$, index $0$, and depth $(s,t)$ is the constant term of an almost meromorphic Jacobi form of index $0$ considered as a polynomial in $\frac{z_2}{\tau_2}$ and $\frac{1}{\tau_2}$.

\subsection{Some modular and elliptic functions}

 For a variable $x$, set $D_x := \frac{1}{\tpi} \frac{\partial}{\partial x}$ and $q_x := e^{2\pi i x}$. Define for $m\in\mathbb{N}:=\{\ell\in \mathbb{Z}: \ell>0\}$ the elliptic functions\footnote{Note that $P_k (w,\tau)$ given in \eqref{eq:Pm} is $P_k (\tpi w,\tau)$ of Section~2.1 of \cite{MTZ}, a $(-\tpi)^{-k}$ multiple of $P_k (q_w,q)$ in Section~3 of  \cite{Zhu}, and a $(\tpi)^{-k}$ multiple of $E_k (w,\tau)$ in Section~2 of \cite{Libgober-Elliptic}.}
\begin{equation}
\label{eq:Pm}
\begin{aligned}
P_{1}(w,\tau):=&-\sum_{n\in\Z\backslash \{0\}}\frac{q_w^n}{1-q^{n}}-\frac{1}{2},
\\
P_{m+1}(w,\tau):=&\frac{\left(-1\right)^m}{m!} D_w^m \left(P_{1}(w,\tau)\right)
=\frac{(-1)^{m+1}}{m!}\sum_{n\in\Z\backslash \{0\}}\frac{n^m q_w^n}{1-q^{n}}.
\end{aligned}
\end{equation}
Note that $P_{m}(w,\tau)={(\tpi w)^{-m}}+O(1)$ in the neighborhood of $w=0$.
Moreover, we require the modular Eisenstein series $G_{k}(\tau)$, defined by
$G_{k}=0$ for $k$ odd whereas for $k\ge 2$ even\footnote{The $G_k$ defined here are precisely the $E_k$ given in Section~2.1 of \cite{MTZ}.}
\begin{align}
 G_{k}(\tau)&=-\frac{ B_{k}}{k!}+\frac{2}{(k-1)!}
\sum\limits_{n\geq 1}\frac{n^{k-1}q^{n}}{1-q^{n}}, \notag 
\end{align}%
where $B_{k}$ is the $k$th Bernoulli number
defined by $(e^z-1)^{-1}
=:\displaystyle{\sum\limits_{k\geq 0}\frac{B_{k}}{k!}z^{k-1}}$. In particular, the first three Bernoulli numbers are given by $B_{0}=1$, $B_{1}=-\frac{1}{2}$, and $B_{2}=\frac{1}{6}$. It is convenient to also define
$G_{0}:=-1$.
Recall that $G_{k}$ is a modular form for $k>2$ and a quasi-modular form for $k=2$. Therefore,
\begin{align}
G_{k}(\gamma \tau )=(c\tau +d)^{k} G_{k}(\tau )-\delta_{k,2} \frac{ c(c\tau +d)}{\tpi},
\notag 
\end{align}
where $\gamma \tau :=\frac{a\tau +b}{c\tau +d}$ for $\gamma =\left(
\begin{smallmatrix}
a & b \\
c & d%
\end{smallmatrix}%
\right) \in \SL_2(\Z)$ and $\delta_{a,b}=1$ if $a=b$ and $0$ otherwise.

These Eisenstein series are related to $P_1$ by
\begin{align}
P_{1}(w,\tau)&=\frac{1}{\tpi w}-\sum_{k\ge 1}G_{k}(\tau)(\tpi w)^{k-1}.
\label{eq:P1Gk}
\end{align}

\medskip
Note that $P_2$ is related to the classical  Weierstrass elliptic function
\cite[Section~2]{Lang}
\[
\wp (w,\tau) := \frac{1}{w} + \sum_{(m,n)\in \mathbb{Z}^2 \setminus \{(0,0)\}} \left(\frac{1}{(w-(m\tau +n))^2} - \frac{1}{(m\tau +n)^2} \right)
\]
 with periods $1$ and $\tau$ by
\begin{align}
P_2(w,\tau)= \frac{1}{(\tpi)^2}\wp(w,\tau)+G_2(\tau).
\notag 
\end{align}
Since $\wp$ is a meromorphic Jacobi form of weight $2$ and index $0$, $P_n$ is a meromorphic Jacobi form of weight $n$ and index $0$ for all $n\ge 3$, while $P_2$ is a quasi-Jacobi form of weight $2$, index $0$, and depth $(0,1)$. That is, for any $\lambda, \mu \in \mathbb{Z}$
\begin{align*}
P_2\left(\gamma.w,\gamma \tau\right)&=(c\tau+ d)^2P_2(w,\tau)- \frac{ c(c\tau +d)}{\tpi}, \\
\notag 
P_2 \left(w +\lambda \tau+\mu, \tau \right) &= P_2 (w,\tau), \notag
\end{align*}
where for $\gamma=\left(\begin{smallmatrix}
a & b\\
c & d
\end{smallmatrix}\right) \in \SL_2(\Z)$ we set $\gamma.w:=\frac{w}{c\tau+d}$.
Lastly, $P_1$ is a quasi-Jacobi form of weight $1$ and index $0$ (and depth $(1,0)$) since for all $\lambda, \mu\in \Z$
\begin{align*}
P_1\left(\gamma.w,\gamma \tau\right)&=(c\tau+ d)P_1(w,\tau)+ \tpi cw, \\
\notag 
P_1(w+\lambda\tau+\mu,\tau)&=P_1(w,\tau)- \lambda .
\notag 
\end{align*}

We also define the elliptic prime form
\begin{align}
K(w,\tau ):=& \exp\left(-P_0(w,\tau)\right)= \tpi w+O\left(w^3\right), \quad \text{where}  \label{eq:Kdef} \\
P_{0}(w,\tau) :=& -\text{Log} (\tpi w)+\sum_{k\ge 1}G_{k}(\tau)\frac{(\tpi w)^{k}}{k}.
\label{eq:P0}
\end{align}
Clearly, $P_1(w,\tau)=D_w \left(\text{Log}\left( K(w,\tau)\right)\right)$ and $P_2(w,\tau)=-D_w^2 \left(\text{Log} (K(w,\tau))\right)$ \cite[Page 34]{Fay}.
The function $K$ is a Jacobi form of weight $-1$ and  index $\half$  with a multiplier system. More precisely, for all $\lambda, \mu\in \Z$
\begin{align*}
K\left(\gamma.w,\gamma \tau\right)&=(c\tau+ d)^{-1}
e^{\frac{\pi i c w^2}{c\tau+ d}}K(w,\tau), \\
\notag 
K(w+\lambda\tau+\mu,\tau)&=(-1)^{\lambda+\mu}q_w^{-\lambda}q^{-\frac{\lambda^2}{2}} K(w,\tau) .
\notag 
\end{align*}
 Additionally, $K$ is expressible in terms of the Jacobi theta function
\begin{align}
\theta_1(w,\tau):=i\sum_{n\in \Z}(-1)^n q^{\half\left(n+\half\right)^2}q_w^{n+\half},
 \label{eq:theta1}
\end{align}
namely
\begin{align*}
K(w,\tau )=\frac{\theta_1(w,\tau)}{D_w \left[\theta_1(w,\tau)\right]_{w=0}}=\frac{\theta_1(w,\tau)}{i \eta^3(\tau)},
\end{align*}%
where $\eta(\tau):=q^{\frac{1}{24}}\prod_{n\ge1}(1-q^n)$ is Dedekind's $\eta$-function, a weight $\half$ modular form.


\subsection{A quasi-Jacobi generating function}
For $w,z\in\C$ and $\tau \in \HH$ define
 \begin{align*}
\widetilde{P}_1(w,z,\tau) :=-\sum_{n\in\Z}\frac{q_w^n}{1-q_z q^n}.
\end{align*}
We note that
$$
\widetilde{P}_1(w,z,\tau)= F_\tau(2\pi iw,2\pi iz),
$$
where $F_\tau$ is given in Section 3 of \cite{Zagier-Periods}.
For $m\in\mathbb{N}_0$ , we also define
 \begin{align}
\widetilde{P}_{m+1}(w,z,\tau) :=\frac{(-1)^{m}}{m!} D_w^m \left(\widetilde{P}_1(w,z,\tau)\right)
=\frac{(-1)^{m+1} }{m!}
\sum_{n\in\Z}\frac{n^m q_w^n}{1-q_zq^n}.
\label{eq:Pmtilde}
\end{align}

The following proposition can be concluded from  \cite{Zagier-Periods}.

\begin{proposition}
\label{prop:P1tilde}
We have
\begin{enumerate}[label=\normalfont(\roman*)]
	\item $\widetilde{P}_1 (w,z,\tau)$ is absolutely convergent for $0<|q|<|q_w|,|q_z|<1$;
	\item $\widetilde{P}_1(w,z,\tau) =\widetilde{P}_1(z,w,\tau) =-\widetilde{P}_1(-w,-z,\tau)$ and $\widetilde{P}_1(w,-w,\tau) =0$;
	\item $\widetilde{P}_1(w,z,\tau) =q_z^{\lambda}\widetilde{P}_1(w+\lambda\tau+\mu,z,\tau)
	=q_w^{\lambda}\widetilde{P}_1(w,z+\lambda\tau+\mu,\tau) $ for all $\lambda,\mu\in \Z$;
	\item $\widetilde{P}_1(w,z,\tau) $ satisfies the differential equation
	\begin{align}
	\left(D_w D_z - D_\tau \right) \widetilde{P}_1(w,z,\tau)=0;
	\notag 
	\end{align}
	\item $\displaystyle{\lim_{z\rightarrow 0} \left (\widetilde{P}_1(w,z,\tau) -\frac{1}{\tpi z}\right)}=P_1(w,\tau)$;
	\item the equality
	       \begin{align}
			\widetilde{P}_1(w,z,\tau) =  \frac{K(w+z,\tau)}{K(w,\tau)K(z,\tau)};
			\label{eq:P1K}
		   \end{align}
	\item for all $\gamma=\left(\begin{smallmatrix}
		a & b\\
		c & d
		\end{smallmatrix}\right)\in \SL_2(\Z)$,
		\begin{align}
		\widetilde{P}_1\left(\gamma.w,\gamma.z,\gamma\tau\right) =
		(c\tau +d)
		e^{\frac{\tpi c wz  }{c\tau+ d}}\,
		\widetilde{P}_1(w,z,\tau).
		\notag 
		\end{align}
\hfill $\Qed$
	\end{enumerate}
\end{proposition}

Proposition~\ref{prop:P1tilde} (v) and \eqref{eq:Pmtilde} imply that
$\widetilde{P}_{m+1}(w,z,\tau)$ has simple poles at $z= \lambda\tau+\mu$ for $\lambda, \mu\in \Z$ with residue $\frac{\lambda^m q_{w}^{-\lambda}}{m!\tpi}$ and no other poles.
It is thus useful to define for $m\in\mathbb{N}_0$
\begin{equation}
\label{eq:PellPm}
P_{m+1, \lambda}\left(w,\tau\right):=
\lim_{z\rightarrow  \lambda\tau+\mu}\left(\widetilde{P}_{m+1}(w,z,\tau)
- \frac{1}{\left(z- \lambda\tau-\mu\right)}\frac{\lambda^m q_{w}^{-\lambda}}{m!\tpi} \right)
=\frac{(-1)^{m+1}}{m!}\sum_{n\in \Z\backslash \{-\lambda\}}\frac{n^mq_w^n}{1-q^{n+\lambda}}.
\end{equation}
We note that $P_{1,\lambda}\left(w,\tau\right)=q_w^{-\lambda}(P_1(w,\tau)+1/2)$ with
\begin{align}
P_{m+1,\lambda}\left(w,\tau\right)&=\frac{(-1)^m}{m!} D_w^m \left(P_{1,\lambda}\left(w,\tau\right)\right).
\notag 
\end{align}
Similarly to \eqref{eq:P1Gk}, we also consider the expansion
\begin{align}
\notag 
P_{1,\lambda}(w,\tau)=:\frac{1}{\tpi w}-\sum_{k\ge 1}G_{k,\lambda}(\tau )(\tpi w)^{k-1},
\end{align}
where we find (see also Section~3 of \cite{Zagier-Periods})
\begin{align}
G_{k,\lambda}(\tau )&=\sum_{j=0}^{k}\frac{\lambda^j}{j!}G_{k-j}(\tau ).
  \label{eq:Gkl}
\end{align}

\medskip
The right side of \eqref{eq:P1K} appears in both \cite{Libgober-Elliptic} and \cite{Oberdieck-Serre} and in other sources as a generating function for quasi-Jacobi forms. Thus \eqref{eq:Kdef} and \eqref{eq:P0} imply that
\begin{align}
\text{Log}\left(\tpi w\widetilde{P}_1(w,z,\tau)\right) =& \text{Log} (\tpi w)+P_0(w,\tau) +P_0(z,\tau)-P_0(z+w,\tau)
\notag
\\
=&  \sum_{k\ge 1}(-1)^{k-1}\left(P_k(z,\tau)-G_k(\tau)\right)\frac{(\tpi w)^k}{k}.
\label{eq:gen1}
\end{align}
Note that $P_k(z,\tau)-G_k(\tau)$ for $k\in\mathbb{N}$ together with $G_2(\tau)$ generate the ring of quasi-Jacobi forms \cite[Proposition 2.9 or 2.10]{Libgober-Elliptic}.
We  also define another generating set $\widetilde{G}_k(z,\tau)$ for $k\ge 1$ together with  $G_2(\tau)$ given by
\cite{Oberdieck-Serre}
\begin{align}
\widetilde{P}_1(w,z,\tau) =: \frac{1}{\tpi w}-\sum_{k\ge 1}\widetilde{G}_k(z,\tau) (\tpi w)^{k-1},
\label{eq:P1Gn}
\end{align}
where we find that for $k\ge 1$,
\begin{equation}
\label{eq:Gktild}
\begin{aligned}
\widetilde{G}_k(z,\tau) =&
-\delta_{k,1}\frac{q_z}{q_z-1}
-\dfrac{B_{k}}{k!}
+\frac{1}{(k-1)!}
\sum_{m,n\ge 1}\left(n^{k-1} q_z^{m}+(-1)^{k}n^{k-1}q_z^{-m} \right)q^{mn},
\notag
\\
=&
-\delta_{k,1}\frac{q_z}{q_z-1}
-\dfrac{B_{k}}{k!}
+\frac{1}{(k-1)!}
\sum_{n\ge 1}\left( \frac{n^{k-1}q_z q^n}{1-q_zq^n}
+(-1)^{k}\frac{n^{k-1}q_z^{-1} q^n}{1-q_z^{-1} q^n} \right).
\end{aligned}
\end{equation}
We note that the function $\widetilde{G}_k (z,\tau)$ equals $-J_k (z,\tau)/k!$, for the functions $J_k$ studied in \cite[display (2)]{Oberdieck-Serre}.
Comparing \eqref{eq:gen1} and \eqref{eq:P1Gn} we  also note that
\begin{align}
\label{eq:G_1tilde}
\widetilde{G}_1(z,\tau)&=-P_1(z,\tau),\\
\label{eq:G_2tilde}
\widetilde{G}_2 (z,\tau) &= \frac{1}{2} \left(P_2 (z,\tau) - P_1 (z,\tau)^2 - G_2 (\tau) \right) .
\end{align}
Defining $ \widetilde{G}_0(z,\tau):=-1$, we find that Proposition~\ref{prop:P1tilde} implies the following (see also \cite{GK-differential} and \cite{Oberdieck-Serre}).
\begin{proposition}
\label{prop:Gtilde}
For all $k\in \mathbb{N}$ and $\lambda, \mu\in \mathbb{Z}$, we have
\begin{enumerate}[label=\normalfont(\roman*)]
		\item $\widetilde{G}_k (z,\tau)$ is absolutely convergent for $0<|q|<|q_z|<1$;
		\item $\widetilde{G}_k(-z,\tau)=(-1)^k\widetilde{G}_k(z,\tau)$;
		\item $\displaystyle{\widetilde{G}_k(z+\lambda\tau+\mu,\tau)
		=\sum_{m=0}^{k}\frac{1}{m!}}(- \lambda)^{m}\widetilde{G}_{k-m}(z,\tau)$;
		\item $\displaystyle{\widetilde{G}_k(\gamma.z,\gamma\tau)
		=(c\tau +d)^k
		\sum_{m=0}^{k}\frac{1}{m!}}\left(\frac{cz}{c\tau+d}\right)^{m}
		\widetilde{G}_{k-m}(z,\tau)$
		for all $\gamma=\left(\begin{smallmatrix}
		a & b\\
		c & d
		\end{smallmatrix}\right)\in \SL_2(\Z)$;
		\item $\widetilde{G}_k (z,\tau)$ has a simple pole at $z=\lambda\tau+\mu$
		with residue $-\frac{1}{\tpi}\frac{(-\lambda)^{k-1}}{(k-1)!}$ and no other poles;
		\item  $\widetilde{G}_k(0,\tau)=G_k(\tau)$ for $k\ge 2$;
		\item $D_\tau\left(\widetilde{G}_k (z,\tau)\right)  = - k D_z \left(\widetilde{G}_{k+1} (z,\tau)\right)$.
\hfill $\Qed$
\end{enumerate}
\end{proposition}

Note that $\tpi w\widetilde{P}_1(w,z,\tau)e^{\tpi w\widetilde{G}_1(z,\tau)}$ is invariant under the action of the Jacobi group. Thus $\text{Log}(\tpi w\widetilde{P}_1(w,z,\tau)e^{\tpi w \widetilde{G}_1(z,\tau)})$ from \eqref{eq:gen1} generates the space of meromorphic Jacobi forms $P_k(z,\tau)-G_k(\tau)$ of index $0$ and weight $k\ge 2$. Alternatively,
\begin{align*}
\tpi w\widetilde{P}_1(w,z,\tau)e^{\tpi w \widetilde{G}_1(z,\tau)}=-\sum_{k\ge 0}K_k (z,\tau)(\tpi w)^{k}
\end{align*}
(see also \eqref{eq:G_1tilde}), where\footnote{The $K_n (z,\tau)$ here equal $-n!K_n (z,\tau)$ for the functions $K_n$ defined in \cite[display (12)]{Oberdieck-Serre}.}
\begin{align}
K_n (z,\tau):=\sum_{m=0}^n\frac{1}{m!}\widetilde{G}_{n-m}(z,\tau) \widetilde{G}_1(z,\tau)^m
\label{eq:Kn}
\end{align}
is a meromorphic Jacobi form of weight $k$ and index $0$
\cite[Proposition $7$]{Oberdieck-Serre}.


\section{Zhu reduction for Jacobi $n$-point functions}\label{Sec:Zhu}

\subsection{Jacobi $n$-point functions}
Let $V$ be a VOA with Virasoro vector $\omega$ of central charge $\mathbf{c}$. Consider $J\in V_{1}$ such that $J(0)$ acts semisimply on $V$.
For  $v_1,\dots,v_n\in V$ and a  weak $V$-module $M$, the Jacobi $n$-point function is
\begin{align}
Z_M^J\left ((v_1,x_1),\dots, (v_n,x_n);z,\tau\right)
:=\tr_{M}Y\left(e^{x_1L(0)}v_1,e^{x_1}\right)\cdots Y\left(e^{x_nL(0)} v_n,e^{x_n}\right)\zeta^{J(0)}q^{L(0)},
\label{eq:npt}
\end{align}
where $\zeta:=q_z=e^{\tpi z}$ as before.
In particular, the Jacobi 1-point function, for $v\in V$, is given by
\begin{align}
Z_M^J\left (v;z,\tau\right)
:=\tr_{M}o(v)\zeta^{J(0)}q^{L(0)}.
\notag 
\end{align}
Define the square bracket operators for $V$ by
\begin{align}
Y[v,x]:=Y\left(e^{xL(0)} v,e^x -1 \right) =: \sum_{n\in \mathbb{Z}} v[n] x^{-n-1}.
\notag 
\end{align}
For $v$ of weight $\wt(v)$ and $k\in \Z$ (see \cite[Lemma~4.3.1]{Zhu}), we have
\begin{align}
\sum_{j\geq 0}\binom{k+\wt(v)-1}{j}v(j)= \sum_{m\ge 0}\frac{k^m}{m!}v[m].
\label{eq:asqround}
\end{align}
The square bracket operators form an isomorphic VOA with Virasoro vector
\begin{align}
\widetilde{\omega}:= \omega-\frac{\mathbf{c}}{24} \vac.
\notag 
\end{align}

\begin{lemma}
\label{lem:Jacnpt}
 The Jacobi $n$-point function $Z_M^J$ obeys the following properties.
\begin{enumerate}[label=(\roman*), font=\normalfont]
\item
We have
\begin{align*}
Z^{J}_{M}\left((v_{1},x_{1}),\cdots ,(v_{n},x_{n});z,\tau \right)
=&Z^{J}_{M}\left(Y[v_{1},x_{1n}]\cdots Y[v_{n-1},x_{(n-1)n}]v_{n};z,\tau \right)
\\
=&Z^{J}_{M}\left(Y[v_{1},x_{1}]\cdots Y[v_{n},x_{n}]\vac;z,\tau \right),
\end{align*}
where $x_{jk} := x_{j}-x_{k}$.
\item[(ii)]  For all adjacent pairs $(j, j+1)$,
\begin{multline*}
x_{j(j+1)}^N Z^{J}_{M}((v_{1},x_{1}),\dots ,(v_{j},x_{j}),(v_{j+1},x_{j+1}),\dots
,(v_{n},x_{n});z,\tau ) \\
=x_{j(j+1)}^N Z^{J}_{M}((v_{1},x_{1}),\dots ,(v_{j+1},x_{j+1}),(v_{j},x_{j}),\dots
,(v_{n},x_{n});z,\tau )
\end{multline*}
for $N\gg 0$.
\item[(iii)] The function $Z^{J}_{M}$ is a function of $x_{jk}$ and is
non-singular at $x_{jk} \neq 0$ for all $j\neq k$.
\item[(iv)] The function $Z^{J}_{M}$ is formally periodic in $x_{j}$ with periods $\tpi$ and $\tpi \tau$ for respective
multipliers $1$ and $e^{\tpi z \alpha_j}$.
\item[(v)] Assume that $J(0)v_j=\alpha_j v_j$ for $\alpha_j \in \mathbb{C}$, $j=1,\dots,n$. Then $Z^{J}_{M}=0$ if $z\sum_{j=1}^n\alpha_j\notin \Z$.
\end{enumerate}
\end{lemma}

\noindent \textbf{Proof.} Part (i) follows from Lemma 1 of \cite{MT} and (ii)--(iv) from Lemma~4 of \cite{MTZ}.
Meanwhile, Part (v) can be deduced by noting that
\begin{align*}
Z_M^J\left ((v_1,x_1),\dots, (v_n,x_n);z,\tau\right)=
\tr_{M}\zeta^{J(0)}Y(e^{x_1 L(0)}v_1,e^{x_1})\cdots Y(e^{x_nL(0)},e^{x_n})q^{L(0)},
\end{align*}
and using that
$$
\zeta^{J(0)}Y(v_j,x_j)\zeta^{-J(0)}=Y(\zeta^{J(0)}v_j,x_j)=e^{\tpi z\alpha_j}Y(v_j,x_j).
$$
\hfill$\Qed$


\subsection{Zhu reduction}
Suppose that $a\in V$ with $L(0)a=\wt(a)a$ and $J(0)a=\alpha a$ for $\alpha \in \mathbb{C}$.
\begin{lemma}
\label{lem:Rec1}
For all $k\in\Z$, we have
\begin{multline*}
\left(1-\zeta^{-\alpha} q^k\right)\trq {M}{a(k+\wt(a)-1)Y\left(e^{x_1L(0)}v_1,e^{x_1}\right)\cdots Y\left(e^{x_nL(0)}v_n, e^{x_n} \right)\zeta^{J(0)}}
\notag
\\
= \sum_{j=1}^{n}
Z_M^J\left ((v_1,x_1),\dots, \left(e^{x_j k}\sum_{m\ge 0}\frac{k^m}{m!}a[m]v_j,x_j\right),\dots, (v_n,x_n);z,\tau\right).
\notag 
\end{multline*}
\end{lemma}

\noindent \textbf{Proof.}  We commute the operator $a(k+\wt(a)-1)$ through the following trace
\begin{multline*}
\trq {M}{a(k+\wt(a)-1)Y\left(e^{x_1 L(0)}v_{1},e^{x_1} \right)\cdots
Y\left(e^{x_n L(0)}v_{n},e^{x_n}\right)\zeta^{J(0)}}\\
\begin{split}
&=\sum_{j=1}^{n}\sum_{r\geq 0}\binom{k+\wt(a)-1}{r}
e^{x_j r} \\
&\quad\times\trq {M}{Y\left(e^{x_1 L(0)}v_{1},e^{x_1}\right)\cdots
Y\left(e^{x_j L(0)}a(r)v_{j},e^{x_j}\right)\cdots
Y\left(e^{x_n L(0)}v_{n},e^{x_n}\right)\zeta^{J(0)}}
\\
&\quad+\zeta^{-\alpha} q^{k}
\trq {M}{a\left(k+\wt(a)-1\right)Y\left(e^{x_1 L(0)}v_{1},e^{x_1}\right)\cdots
Y\left(e^{x_n L(0)}v_{n},e^{x_n}\right)\zeta^{J(0)}}.\hspace{1.8cm}
\end{split}
\end{multline*}%
Using \eqref{eq:asqround}  the result follows.
\hfill $\Qed$
\medskip

Lemma~\ref{lem:Rec1} immediately implies the following corollary.
\begin{corollary}
\label{cor:ZeroRes}
Let $J(0)a=\alpha a$.
If $\alpha z=\lambda\tau+\mu
\in {\mathbb{Z}\tau}+\Z$, then
\begin{align}
\sum_{j=1}^{n}Z_M^J\left ((v_1,x_1),\dots, \left(e^{x_j \lambda}\sum_{m\ge 0}
\frac{\lambda^m}{m!}a[m]v_j,x_j\right),\dots, (v_n,x_n);z,\tau\right)=0.
\label{eq:ZeroRes}
\end{align}
\end{corollary}

We can now state the first Zhu reduction formula for formal Jacobi $n$-point functions.
\begin{proposition}
\label{prop:Zhured}
 Let $a,v_{1},\dots ,v_{n}\in V$ with $J(0)a=\alpha a$, $\alpha \in \mathbb{C}$.  For $\alpha z\notin {\mathbb{Z}\tau} +\mathbb{Z}$, we have
\begin{align}
&Z_M^J\left ((a,y),(v_1,x_1),\dots, (v_n,x_n);z,\tau\right)\notag
\\
&=
 \sum_{j=1}^{n}\sum_{m\ge 0}\widetilde{P}_{m+1} \left(\frac{y-x_j}{\tpi} ,\alpha z,\tau \right)
Z_M^J\left ((v_1,x_1),\dots, (a[m]v_j,x_j),\dots, (v_n,x_n);z,\tau\right).
\label{eq:ZhuRed}
\end{align}
\end{proposition}

\noindent \textbf{Proof.}
From Lemma~\ref{lem:Rec1}, we find
\begin{align*}
&\hspace{-17mm}Z_M^J\left ((a,y),(v_1,x_1),\dots, (v_n,x_n);z,\tau\right)\\
&=\sum_{k\in \Z}e^{-yk}\trq {M}{a(k+\wt(a)-1)Y(e^{x_1 L(0)}v_1,e^{x_1})\cdots Y(e^{x_n L(0)}v_n,e^{x_n})\zeta^{J(0)}}\\
&= \sum_{j=1}^{n}
\sum_{m\ge 0}H_m(e^y,e^{x_j},\zeta)Z_M^J\left (\left(v_1 ,x_1\right) ,\dots, \left(a[m]v_j,x_j\right),\dots ,\left(v_n ,x_n\right);z,\tau\right),
\end{align*}
where
\begin{align*}
H_m \left(e^y,e^{x_j},\zeta \right):=&\frac{1}{m!}\sum_{k\in \Z}\frac{k^m e^{k\left(x_j -y\right)}}{1-\zeta^{-\alpha} q^k}=
\widetilde{P}_{m+1} \left(\frac{y-x_j}{\tpi} ,\alpha z,\tau \right),
\end{align*}
using that $\widetilde{P}_{m+1}(-w,-z,\tau) =(-1)^{m+1}\widetilde{P}_{m+1}(w,z,\tau)
$.
\hfill $\Qed$
\medskip

The $\widetilde{P}_{m+1}$-terms in \eqref{eq:ZhuRed} have simple poles at $\alpha z\in {\mathbb{Z}\tau}+\Z$ but thanks to Corollary~\ref{cor:ZeroRes} the residue at each pole is zero as follows. Consider the principal part at $\alpha z =\lambda\tau+\mu$ defined by
\begin{align}\notag %
\PP\left(\widetilde{P}_{m+1} \left(\frac{y-x_j}{\tpi},\alpha z,\tau \right)\right)
:=\widetilde{P}_{m+1}\left(\frac{y-x_j}{\tpi},\alpha z,\tau \right)
-\frac{\lambda^m e^{\lambda \left(x_j -y\right)}}{m!\tpi (\alpha z-\lambda\tau-\mu)}.
\end{align}
The right side of \eqref{eq:ZhuRed} can be written as
\begin{align*}
&\hspace{-2mm}\sum_{j=1}^{n}\sum_{m\ge 0}
\PP\left(\widetilde{P}_{m+1} \left(\frac{y-x_j}{\tpi},\alpha z,\tau \right)\right)
Z_M^J\left (\left(v_1 ,x_1\right), \dots, (a[m]v_j,x_j),\dots ,\left(v_n ,x_n\right);z,\tau\right)\\
&\quad +\frac{e^{-y\lambda}}{\tpi (\alpha z-\lambda\tau-\mu)}
\sum_{j=1}^{n}\sum_{m\ge 0}
\frac{\lambda^m e^{x_j \lambda}}{m!}
Z_M^J\left (\left(v_1 ,x_1\right), \dots, (a[m]v_j,x_j),\dots ,\left(v_n ,x_n\right);z,\tau\right)\\
&\hspace{-4mm}=
\sum_{j=1}^{n}\sum_{m\ge 0}
\PP\left(\widetilde{P}_{m+1} \left(\frac{y-x_j}{\tpi},\alpha z,\tau \right)\right)
Z_M^J\left (\left(v_1 ,x_1\right), \dots, (a[m]v_j,x_j),\dots ,\left(v_n ,x_n\right);z,\tau\right)\\
&\hspace{-4mm}+
\frac{\left(1-\zeta^{-\alpha} q^\lambda \right)}{\tpi (\alpha z-\lambda\tau-\mu)}
e^{-y\lambda}\trq {M}{a(\lambda+\wt(a)-1)Y\left(e^{x_1 L(0)}v_1,e^{x_1}\right)\cdots Y\left(e^{x_n L(0)}v_n,e^{x_n}\right)\zeta^{J(0)}},
\end{align*}
by Lemma~\ref{lem:Rec1}. This residue at $\alpha z = \lambda\tau+\mu$ is zero and we establish the following result.
\begin{proposition}
\label{prop:Zhured0}
 Let $a,v_{1},\dots ,v_{n}\in V$ with $J(0)a=\alpha a$.  For $\alpha z=\lambda\tau+\mu\in {\mathbb{Z}\tau+\Z}$, we have
\begin{align}
&Z_M^J\left((a,y),(v_1,x_1),\dots, (v_n,x_n);z,\tau\right)\notag
\\
&\quad =
e^{-y\lambda}\tr_{M}a(\lambda+\wt(a)-1)Y(e^{x_1 L(0)}v_1,e^{x_1})\cdots Y(e^{x_n L(0)}v_n,e^{x_n})
\zeta^{J(0)}q^{L(0)}
\notag
\\
&\quad \quad +\sum_{j=1}^{n}\sum_{m\ge 0}P_{m+1,\lambda}\left( \frac{y-x_j}{\tpi} ,\tau\right)
Z_M^J \left(\left(v_1 ,x_1\right), \dots, (a[m]v_j,x_j),\dots ,\left(v_n ,x_n\right);z,\tau \right),
\label{eq:ZhuRed0}
\end{align}
with  $P_{m+1,\lambda}\left(w ,\tau\right)$ defined in \eqref{eq:PellPm}.
\hfill $\Qed$
\end{proposition}

We are in position to describe the second Zhu reduction formula for Jacobi $n$-point functions.
\begin{proposition}
\label{prop:apnpt}
 Let $a,v_{1},\dots v_{n}\in V$ with $J(0)a=\alpha a$. For $N\ge 1$ and $\alpha z \notin  {\mathbb{Z}\tau}+\Z $, we have
\begin{align}
&Z_M^J\left ((a[-N]v_1,x_1),\dots, (v_n,x_n);z,\tau\right)
\notag
\\
&\quad =
\sum_{m\ge 0}(-1)^{m+1}\binom{m+N-1}{m}\widetilde{G}_{m+N}(\alpha z,\tau)
Z_M^J\left ((a[m]v_1,x_1),\dots, (v_n,x_n);z,\tau\right)\notag
\\
&\quad \quad+
\sum_{j=2}^{n}\sum_{m\ge 0}
(-1)^{N+1}
\binom{m+N-1}{m}
\widetilde{P}_{m+N} \left(\frac{x_1-x_j}{\tpi} ,\alpha z,\tau \right)\notag
\\
&\quad \qquad\qquad\times
Z_M^J\left ((v_1,x_1),\dots, (a[m]v_j,x_j),\dots, (v_n,x_n);z,\tau\right).
\label{eq:2ZhuRed}
\end{align}

\end{proposition}

\noindent \textbf{Proof.}
Using (i) of Lemma~\ref{lem:Jacnpt} and the associativity of VOAs, we find that
\begin{equation}\label{FY[v,z]}
Z_{M}^J((Y[v,y]v_{1},x_{1}),\dots ,(v_{n},x_{n});z,\tau )
=Z_{M}^J((v,x_{1}+y),(v_{1},x_{1}),\dots ,(v_{n},x_{n});z,\tau ).
\end{equation}%
Expanding the left side of \eqref{FY[v,z]} in $y$ gives that the coefficient of
$y^{N-1}$ equals
\[
Z_{M}^J \left( \left(v[-N]v_{1},x_{1} \right), \left(v_2 ,x_2 \right), \dots , \left(v_{n},x_{n}\right);z,\tau \right).
\]
We can
compare this to the expansion of $y^{N-1}$ in the right side of Proposition~\ref{prop:Zhured}.
From \eqref{eq:ZhuRed} we see that
the coefficient of $y^{N-1}$ in
$\widetilde{P}_{m+1}((x_1+y-x_1)/\tpi ,\alpha z,\tau)$
is $(-1)^{m+1}\binom{m+N-1}{m}\widetilde{G}_{m+N}(\alpha z,\tau )$, using \eqref{eq:P1Gn} and that for
$\widetilde{P}_{m+1}((y+x_1-x_j)/\tpi ,\alpha z,\tau)$
is $(-1)^{p+1}\binom{m+N-1}{m}\widetilde{P}_{m+N}((x_1-x_j)/\tpi ,\alpha z,\tau )$
for $j\neq 1$. Thus \eqref{eq:2ZhuRed} follows.
\hfill $\Qed$
\medskip

Propositions~\ref{prop:Zhured0} and \ref{prop:apnpt} imply the next result.
\begin{proposition}
\label{prop:apnpt0}
 Let $a,v_{1},\dots v_{n}\in V$ with $J(0)a=\alpha a$. For $N\geq 1$ and  $\alpha z = \lambda\tau+\mu\in {\mathbb{Z}\tau}+\Z$, we have
\begin{align}
&Z_M^J\left ((a[-N]v_1,x_1),\dots, (v_n,x_n);z,\tau\right)\notag
\\
&\quad =
(-1)^{N+1}\frac{\lambda^{N-1}}{(N-1)!}
\tr _{M}a(\lambda+\wt(a)-1)Y\left(e^{x_1 L(0)}v_1,e^{x_1}\right)\cdots Y\left(e^{x_n L(0)}v_n,e^{x_n}\right)
\zeta^{J(0)}q^{L(0)}
\notag
\\
&\quad \quad +\sum_{m\ge 0}(-1)^{m+1}\binom{m+N-1}{m}
{G}_{m+N,\lambda}(\tau)
Z_M^J\left ((a[m]v_1,x_1),\dots, (v_n,x_n);z,\tau\right)\notag
\\
&\quad \quad+
\sum_{j=2}^{n}\sum_{m\ge 0}
(-1)^{N+1}
\binom{m+N-1}{m}
P_{m+N,\lambda}\left(\frac{x_1-x_{j}}{\tpi},\tau \right)\notag
\\
&\quad \qquad\qquad\times
Z_M^J\left ((v_1,x_1),\dots, (a[m]v_j,x_j),\dots, (v_n,x_n);z,\tau\right),
\notag 
\end{align}
for ${G}_{k, \lambda}$ given in \eqref{eq:Gkl}.
\hfill $\Qed$
\end{proposition}

\begin{remark}
In the case $\alpha=0$ we have that $\lambda=\mu=0$ and Propositions~\ref{prop:Zhured0} and \ref{prop:apnpt0} imply the standard results of \cite{Zhu} or \cite{MTZ} with $a(\lambda+\wt(a)-1)=o(a)$.
\end{remark}

Theorem \ref{thm:intro1} now follows as a corollary to Propositions \ref{prop:Zhured}--\ref{prop:apnpt0}.


\subsection{Zhu reduction with a shifted Virasoro vector}
In this section, we show that Corollary~\ref{cor:ZeroRes} and
Proposition~\ref{prop:Zhured0} are related to previously known results based on an appropriate shifted Virasoro vector. Suppose that $J(0)a=\alpha a$ for $\alpha\not \in \mathbb{Z}\setminus \{0\}$, and define $g\in\Aut(V)$ by
\begin{align}
g:=e^{\frac{\tpi\mu}{\alpha}J(0)},
\notag 
\end{align}
for $\mu\in \Z$ for which $ga=a$.
Then Corollary~\ref{cor:ZeroRes} follows from
 Proposition~6 of \cite{MTZ} which states that
\begin{align}
\sum_{j=1}^{n}\tr_{M}Y\left(e^{x_1 L(0)}v_1,e^{x_1}\right)\cdots
Y\left(e^{x_j L(0)}a[0]v_j,e^{x_j}\right)\cdots
 Y \left(e^{x_n L(0)}v_n,e^{x_n} \right)\,gq^{L(0)}=0.
\label{eq:Zg0}
\end{align}
Consider the shifted Virasoro vector
\begin{align}
\omega_{h}:=\omega +h(-2)\vac,
\notag 
\end{align}
where  $h=-\frac{\lambda}{\alpha}J$ for $\lambda\in \Z$.
The shifted grading operator is
\begin{align*}
L_{h}(0):=L(0)-h(0)=L(0)+\frac{\lambda}{\alpha}J(0).
\end{align*}
 Denote the square bracket vertex operator for the shifted Virasoro vector by
\begin{align}
Y[v,x]_h:=Y\left(e^{xL_h(0)}v,e^x -1\right) =: \sum_{n\in \Z} v[n]_hx^{-n-1}.
\notag 
\end{align}
Hence $Y[a,x]_h=e^{x\lambda} Y[a,x]$,
or equivalently,
\begin{align}
a[n]_h=\sum_{m\ge 0}\frac{\lambda^m}{m!}a[n+m].
\label{eq:anh}
\end{align}
Next consider \eqref{eq:Zg0} with the shifted Virasoro grading.
With $J(0)v_j=\alpha_j v_j$ for $j=1,\dots,n$ we find
\begin{align*}
0
=&
\sum_{j=1}^{n}\tr_{M}Y\left(e^{x_1 L_{h}(0)}v_1,e^{x_1}\right)\cdots
Y\left(e^{x_j L_{h}(0)}a[0]_hv_j,e^{x_j}\right)\cdots
 Y\left(e^{x_n L_{h}(0)}v_n,e^{x_n}\right)\,gq^{L_h(0)}\\
=&
\prod_{r=1}^{n}e^{x_r \frac{\lambda\alpha_r}{\alpha}}
\sum_{j=1}^{n}\tr_{M}Y\left(e^{x_1 L(0)}v_1,e^{x_1}\right)\cdots
Y\left(e^{x_j \lambda+L(0)}a[0]_hv_j,e^{x_j}\right)\cdots
 Y\left(e^{x_n L(0)}v_n,e^{x_n}\right)\,gq^{L_h(0)}.
\end{align*}
Thus, using \eqref{eq:anh} for $n=0$, we find that
\begin{align*}
0&=\sum_{j=1}^{n}\tr_{M}Y\left(e^{x_1 L(0)}v_1,e^{x_1}\right)\cdots
Y\left(e^{x_j \lambda+L(0)}
\sum_{m\ge 0}\frac{\lambda^m}{m!}a[m]v_j,e^{x_j}\right)\cdots
 Y\left(e^{x_n L(0)}v_n,e^{x_n}\right)\,gq^{L_h(0)}\\
&=
\sum_{j=1}^{n}Z_M^J\left ((v_1,x_1),\dots, \left(e^{x_j \lambda}
\sum_{m\ge 0}\frac{\lambda^m}{m!}a[m]v_j,x_j\right),\dots, (v_n,x_n);\frac{1}{\alpha}(\lambda\tau+\mu),\tau\right),
\end{align*}
i.e., we recover \eqref{eq:ZeroRes} of Corollary~\ref{cor:ZeroRes}.

\medskip
In a similar fashion, we can relate Proposition~\ref{prop:Zhured0} to Theorem~2 of \cite{MTZ} for the above shifted Virasoro grading $L_h(0)$ and with $g=e^{\frac{\tpi\mu}{\alpha}J(0)}$.  Theorem~2 of \cite{MTZ} states that
\begin{align}
&\tr_{M}Y \left(e^{yL_{h}(0)}a,e^y \right)Y \left(e^{x_1 L_{h}(0)}v_1,e^{x_1}\right)\cdots
 Y\left(e^{x_n L_{h}(0)}v_n,e^{x_n}\right)\,gq^{L_{h}(0)}
\notag
\\
&=
\tr_{M}o_{h}(a)Y\left(e^{x_1 L_{h}(0)}v_1,e^{x_1}\right)\cdots
 Y\left(e^{x_n L_{h}(0)}v_n,e^{x_n}\right)\,gq^{L_{h}(0)}
\notag
\\
&\quad +\sum_{j=1}^{n}\sum_{r\ge 0}P_{r+1} \left(\frac{y-x_j}{\tpi} ,\tau \right)
\notag
\\
&  \qquad
\times\tr_{M}
Y\left(e^{x_1 L_{h}(0)}v_1,e^{x_1}\right)\cdots
Y\left(e^{x_j L_{h}(0)}a[r]_{h}v_j,e^{x_j}\right)\cdots
 Y\left(e^{x_n L_{h}(0)}v_n,e^{x_n}\right)\,gq^{L_{h}(0)},
\label{eq:Zhuorb}
\end{align}
where $o_{h}(a):=a(\wt_{h}(a)-1)=a(\mu+\wt(a)-1)$. With $J(0)v_j=\alpha_j v_j$ for $j=1,\dots,n$ the left side of \eqref{eq:Zhuorb} equals
\begin{align}
e^{y \lambda} \prod_{j=1}^{n}e^{x_j \frac{\lambda \alpha_j}{\alpha}}
Z_M^J\left ((a,y),(v_1,x_1),\dots, (v_n,x_n);\frac{1}{\alpha}\left(\lambda \tau + \mu \right),\tau\right),
\label{eq:LHS}
\end{align}
whereas the right side of \eqref{eq:Zhuorb} gives
\begin{equation}
\label{eq:RHS}
\begin{aligned}
&\prod_{j=1}^{n}e^{x_j \frac{\lambda \alpha_j}{\alpha}}
\tr_{M}a(\lambda +\wt(a)-1)Y\left(e^{x_1 L(0)}v_1,e^{x_1}\right)\cdots
 Y\left(e^{x_n L_{h}(0)}v_n,e^{x_n}\right)e^{\tpi\frac{1}{\alpha}(\mu+\lambda\tau)J(0)}q^{L(0)}
\\
&+\prod_{r=1}^{n}e^{x_r\frac{\lambda \alpha_r}{\alpha}}
\sum_{j=1}^{n}\sum_{r\ge 0}P_{r+1}\left(\frac{y-x_j}{\tpi} ,\tau \right) \\
 &\qquad \times Z_M^J\left ((v_1,x_1),\dots, \left(e^{x_j \lambda}a[r]_{h}v_j,x_j \right),\dots, (v_n,x_n);\frac{1}{\alpha}\left(\lambda \tau +\mu \right),\tau\right).
\end{aligned}
\end{equation}
Next we note that
\begin{align*}
\sum_{r\ge 0}P_{r+1}\left(\frac{w}{\tpi} ,\tau \right)a[r]_{h}=& - \sum_{r,s\ge 0}\sum_{k\in\Z\backslash \{0\}}\frac{(-k)^r \lambda^s}{r!s!}\frac{ e^{wk}}{1-q^{k}}a[r+s]\\
=& \sum_{m\ge 0} \frac{(-1)^{m+1}}{m!}\sum_{k\in\Z\backslash \{0\}}\frac{(k- \lambda)^m e^{wk}}{1-q^{k}}a[m]
= \sum_{m\ge 0}P_{m+1,\lambda}\left(\frac{w}{\tpi},\tau\right) e^{w \lambda} a[m].
\end{align*}
Hence the identity \eqref{eq:ZhuRed0} follows from dividing \eqref{eq:LHS} and  \eqref{eq:RHS} by $e^{y\lambda} \prod_{j=1}^{n}e^{x_j\frac{\lambda \alpha_j}{\alpha}}$.


\section{Differential operators on Jacobi forms}\label{Sec:DiffOpsOnJacobi}
In this section we consider a generalization of  differential operators on $J_{k,m}$, the space  of Jacobi forms of weight $k$ and index $m$, introduced in \cite{Oberdieck-Serre}. We investigate how these operators appear in a number of vertex operator constructions for Jacobi $n$-point functions in the subsequent sections.

For $\alpha\in\Z\setminus  \{0\}$, $m,k\in \N$ and $A,B\in \C$, with $A\neq 4mB$,   define the differential operator
\begin{align}
\notag
 \mathcal{M}=&\mathcal{M}_{k, \alpha} = \mathcal{M}_{(A,B),k, \alpha ,m} \\
\label{eq:Moperator}
:= & {\vartheta}_k + \frac{1}{A-4mB}\left[ B \left(D_z^2 + 2m G_2 (\tau)\right) + \frac{A}{\alpha}\left( \widetilde{G}_1 (\alpha z,\tau) D_z - \frac{2m}{\alpha} \widetilde{G}_2 (\alpha z,\tau) \right)\right],
\end{align}
where ${\vartheta}={\vartheta}_k := D_\tau +kG_2 (\tau)$ is the Serre modular derivative which maps modular forms of weight $k$ to modular forms of weight $k+2$.
(We often use the notation $\vartheta$ without subscript if it is applied to forms with the weight not specified or to functions which are linear combinations of forms with different weights.)
We also define operators for particular values of $(A,B)$ as  follows:
 \begin{equation}
 \label{eq:operators}
 \mathcal{H}=\mathcal{H}_k := \mathcal{M}_{(0,1),k, \alpha ,m}, \quad \mathcal{T}=\mathcal{T}_{k,\alpha} := \mathcal{M}_{(1,0),k, \alpha ,m}, \quad \operatorname{and} \quad \mathcal{S}= \mathcal{S}_{k,\alpha} := \mathcal{M}_{(1,1),k, \alpha ,m}.
 \end{equation}
We remark that $ \mathcal{M},\mathcal{H},\mathcal{T},\mathcal{S}$ are linearly dependent with
\begin{align*}
\mathcal{M}=\frac{1}{A-4mB}(A\mathcal{T}-4mB\mathcal{H}),\quad
\mathcal{T}=4m\mathcal{H}+(1-4m)\mathcal{S}.
\end{align*}
In particular,
\begin{equation}
 \label{eq:Heat}
\mathcal{H}_k={\vartheta}_k - \frac{1}{4m}D_z^2 -\frac{1}{2} G_2 (\tau)=D_\tau- \frac{1}{4m}D_z^2+\left(k-\half\right)G_2(\tau),
\end{equation}
is the well-known modified heat operator which maps (weak) Jacobi forms of weight $k$ and index $m$ to (weak) Jacobi forms of weight $k+2$ and the same index   \cite{EZ}.
Furthermore,
\begin{equation}
 \label{eq:Top}
\mathcal{T}_{k,\alpha}={\vartheta}_k+\frac{1}{\alpha}\widetilde{G}_1 (\alpha z,\tau) D_z - \frac{2m}{\alpha^2} \widetilde{G}_2 (\alpha z,\tau).
\end{equation}
For $\alpha= \pm 1$ we  find
\begin{align*}
\mathcal{S}_{2k,\pm 1}=\partial^J:={\vartheta}_k + \frac{1}{1-4m} \left[D_z^2 +\widetilde{G}_1 ( z,\tau) D_z + 2m \left(G_2 (\tau) - \widetilde{G}_2 (z,\tau)\right) \right],
\end{align*}
the generalized Serre derivative $\partial^J$ for even weight $2k$ index $m$ Jacobi forms introduced in \cite{Oberdieck-Serre}.
Since $\mathcal{H}:J_{k,m, \chi}\rightarrow J_{k+2,m,\chi}$ and $\partial^J:J_{2k,m,\chi}\rightarrow J_{2k+2,m,\chi}$, it is natural to consider the action of the general differential operator $\mathcal{M}$ of \eqref{eq:Moperator} on Jacobi forms.

\begin{lemma}\label{lem:LemPreserve}
  Let $\alpha \in \mathbb{Z}$.
 \begin{enumerate}[label=\normalfont(\roman*),itemsep=-1ex]
	\item The operator $\mathcal{M}_{k,\alpha}$ maps forms transforming like \eqref{eq:Jactr} of weight $k$ with multiplier $\chi$ to forms of weight $k+2$ with  multiplier $\chi$.
  \item
	Assume that $\frac{2}{\alpha N_2}\in\Z$, $\chi\left(\frac{2}{\alpha}, 0\right)=e^{\tpi \frac{a}{N}}$ for $a,N\in \N$ with $\gcd(a,N)=1$ and $N$ odd, and $\chi\left(\begin{smallmatrix}-1&0\\0&-1\end{smallmatrix}\right)=(-1)^k$. Then we have
   \begin{align}
   \label{eq:OnJacobiForms}
   \mathcal{M}_{k,\alpha} \colon J_{k,m, \chi} \to J_{k+2,m, \chi}.
   \end{align}
 \end{enumerate}
\end{lemma}

\noindent \textbf{Proof.}
It is well-known that $\mathcal{H}_{k}$ satisfies the properties (i) and (ii). Hence since $\mathcal{M}_{k,\alpha}$ is a linear combination of $\mathcal{H}_{k}$ and $\mathcal{T}_{k,\alpha}$,  it suffices to prove the results for $\mathcal{T}_{k,\alpha}$.
Let $\phi ( z,  \tau) $ be transforming as  \eqref{eq:Jactr}, i.e., with no condition on holomorphicity.
Note that
 \[
 D_\tau \left(\phi (\gamma z, \gamma \tau)\right) =\frac{1}{(c\tau +d)^2}\left[D_\tau\left( \phi (z,\tau)\right)\right]_{\substack{\tau =\gamma \tau \\ z=\gamma z}} -\frac{c z}{(c\tau +d)^2}\left[D_z \left(\phi (z,\tau ) \right)\right]_{\substack{\tau =\gamma \tau \\ z=\gamma z}} ,
 \]
 while in the variable $z$ we have $D_z \left( \phi (\gamma z, \gamma \tau)\right) = \frac{1}{c\tau +d}\left[D_z\left( \phi ( z,\tau)\right)\right]_{\substack{\tau =\gamma \tau \\ z=\gamma z}}$.
Then a straightforward calculation using Proposition~\ref{prop:Gtilde}~(iv) reveals that
 \begin{align*}
 & \mathcal{T}_{k,\alpha} \left(\phi \vert_{k,m} \gamma  \right)
=\mathcal{T}_{k,\alpha} \left(\phi \right) \vert_{k+2,m} \gamma .
 \end{align*}
We also have the total derivatives
 \[
 D_\tau\left( \phi (z+\tau ,\tau)\right) =\left[D_\tau \left(\phi (z,\tau)\right)\right]_{z=z+\tau} + \left[D_z \left(\phi (z,\tau)\right)\right]_{z=z+\tau}
 \]
 and $D_z \left(\phi (z+\tau ,\tau)\right) = \left[D_z \left(\phi (z,\tau)\right)\right]_{z=z+\tau}$ so that
using Proposition~\ref{prop:Gtilde} (iii) we find
 \begin{align*}
 &\mathcal{T}_{k,\alpha} \left(e^{2\pi i m(\tau +2z)} \phi (z+\tau ,\tau)\right)
 =e^{2\pi i m(\tau +2z)} \left[ \mathcal{T}_{k,\alpha} \left(\phi (z,\tau) \right)\right]_{z=z+\tau} ,
 \end{align*}
  which proves  statement (i). It is easy to see that a weight $k$ Jacobi-like form equipped with a multiplier system is mapped to a weight $k+2$ Jacobi-like form with the same multiplier.

  For (ii), assume $\phi \in J_{k,m, \chi}$
	and $\alpha$ satisfy the conditions of the lemma. By (i) it suffices to determine whether $\mathcal{T}_{k,\alpha}$
	introduces poles via the $\widetilde{G}_1 (\alpha z,\tau)$ and $\widetilde{G}_2 (\alpha z,\tau)$ terms. We shift $z\mapsto z+\frac{1}{\alpha} \left(\lambda \tau +\mu \right)$, for $\lambda ,\mu \in \mathbb{Z}$, in
 \[
 \left(\widetilde{G}_1 (\alpha z,\tau)D_z - \frac{2m}{\alpha} \widetilde{G}_2 (\alpha z,\tau)\right) \phi (z,\tau)
 \]
 to get
 \[
 \left(\widetilde{G}_1 \left(\alpha z +\lambda \tau +\mu ,\tau \right)D_z - \frac{2m}{\alpha} \widetilde{G}_2 \left(z + \lambda \tau +\mu ,\tau \right)\right) \phi \left(z +\frac{\lambda}{\alpha}\tau +\frac{\mu}{\alpha},\tau \right).
 \]
 We are now interested in whether this expression has poles at $z=0$. Using Proposition \ref{prop:Gtilde} (iv), we find that
 \begin{align*}
 &\left(\widetilde{G}_1 \left(\alpha z +\lambda \tau +\mu ,\tau \right)D_z - \frac{2m}{\alpha} \widetilde{G}_2 \left(z + \lambda \tau +\mu ,\tau \right)  \right) \phi \left(z +\frac{\lambda}{\alpha}\tau +\frac{\mu}{\alpha},\tau \right) \\
 &\quad = \left( \lambda D_z +\widetilde{G}_1 (\alpha z,\tau)D_z + \frac{m\lambda^2}{\alpha} + \frac{2m\lambda}{\alpha} \widetilde{G}_1 (\alpha z,\tau) - \frac{2m}{\alpha} \widetilde{G}_2 (\alpha z, \tau) \right) \phi \left(z +\frac{\lambda}{\alpha}\tau +\frac{\mu}{\alpha},\tau \right).
 \end{align*}
 Since $\widetilde{G}_2 (\alpha z, \tau)$ does not have a pole at $z=0$, we are only concerned with the term
 \begin{align*}
 &\widetilde{G}_1 (\alpha z,\tau) \left(D_z + \frac{2m\lambda}{\alpha} \right) \phi \left(z +\frac{\lambda}{\alpha}\tau +\frac{\mu}{\alpha},\tau \right).
 \end{align*}
 Using the Fourier expansion of $\phi(z,\tau)$, we find that
 \begin{align}
 \left[\left(D_z + \frac{2m\lambda}{\alpha} \right) \phi \left(z +\frac{\lambda}{\alpha}\tau +\frac{\mu}{\alpha},\tau \right)\right]_{z=0}
 &= \left[\left(D_z + \frac{2m\lambda}{\alpha} \right)
\sum_{n\in \N_0+\rho_{1}} \sum_{r\in\Z+\rho_{2}\atop{r^2\leq 4nm}}
c(n,r)e^{2\pi i \frac{\mu r}{\alpha}} q^{n+ \frac{\lambda r}{\alpha}} \zeta^{r} \right]_{z=0} \notag \\
 &= \sum_{n\in \N_0+\rho_{1}}
 \sum_{r\in\Z+\rho_{2}\atop{r^2\leq 4nm}} \left(r+\frac{2m\lambda}{\alpha}\right)
c(n,r)e^{2\pi i \frac{\mu r}{\alpha}} q^{n+ \frac{\lambda r}{\alpha}}, \label{eq:PoleProof1}
 \end{align}
where $\rho_{j}\equiv\frac{a_j}{N_j}\pmod \Z$ with $0\le \rho_j<1$. To finish the proof, we show that \eqref{eq:PoleProof1} vanishes under the conditions of the lemma.

The transformation law of $\phi$ and properties of $\chi$ imply that for $\ell\in\Z$  
\[
c(n,-r)=(-1)^k\chi\left(\begin{matrix} -1&0\\0&-1\end{matrix}\right)c(n,r),
\quad
c\left(n+\ell r+m\ell^2,r+2m\ell \right)= \chi(-\ell,0)c(n,r).
\]
Making the change of variables $n\mapsto n+\frac{2\lambda r}{\alpha}+\frac{4m\lambda^2}{\alpha^2}\in \N_{0}+\rho_{1}$ and  $r\mapsto -r-\frac{4m\lambda}{\alpha}\in\Z+\rho_{2}$  (using $\frac{2}{\alpha N_2}\in\Z$), \eqref{eq:PoleProof1} equals 
\begin{align*}
&\sum_{n,r} \left(-r-\frac{2m\lambda}{\alpha}\right) c\left(n+\frac{2\lambda r}{\alpha}+\frac{4m\lambda ^2}{\alpha^2},-r-\frac{4m\lambda}{\alpha}\right)e^{2\pi i\frac{\mu}{\alpha}\left(-r-\frac{4m\lambda}{\alpha}\right)}q^{n+\frac{2\lambda r}{\alpha}+\frac{4m\lambda^2}{\alpha^2}+\frac{\lambda}{\alpha}\left(-r-\frac{4m\lambda}{\alpha}\right)}\\
& = (-1)^{k+1}\chi\left(\begin{matrix} -1&0\\0&-1\end{matrix}\right)\sum_{n,r}\left(r+\frac{2m\lambda}{\alpha}\right) c\left(n+\frac{2\lambda r}{\alpha}+\frac{4m\lambda^2}{\alpha^2},r+\frac{4m\lambda}{\alpha}\right)e^{-\frac{2\pi i\mu r}{\alpha}-8\pi i \frac{\mu m\lambda}{\alpha^2}} q^{n+\frac{\lambda r}{\alpha}}.
\end{align*}	
Since $\frac{2\lambda}{\alpha},\frac{2\mu r}{\alpha} \in \Z$, we obtain
\[c\left(n+\frac{2\lambda r}{\alpha} +\frac{4m\lambda^2}{\alpha^2},r+\frac{4m\lambda}{\alpha}\right) = e^{-\tpi \frac{a\lambda}{N}} c(n,r),
\]
using $\chi\left(\frac{2}{\alpha}, 0\right)=e^{\tpi \frac{a}{N}}$,
and
\[ e^{-\frac{8\pi i\mu m\lambda}{\alpha^2}}=1 ,\quad  e^{\frac{2\pi i\mu r}{\alpha}}=(-1)^{\frac{2\mu r}{\alpha}}. \]
Thus
\begin{align*}
\quad&\sum_{n,r}\left(r+\frac{2m\lambda}{\alpha}\right) c(n,r) (-1)^{\frac{2\mu r}{\alpha}} q^{n+\frac{\lambda r}{\alpha}}\\
&\qquad \quad\quad=(-1)^{k+1}\chi\left(\begin{matrix} -1&0\\0&-1\end{matrix}\right)
e^{-\tpi \frac{a\lambda}{N}}
 \sum_{n,r} \left(r+\frac{2m\lambda}{\alpha}\right) c(n,r) (-1)^{\frac{2\mu r}{\alpha}} q^{n+\frac{\lambda r}{\alpha}}.
\end{align*}
Hence, provided that
\begin{eqnarray}
\chi\left(\begin{matrix} -1&0\\0&-1\end{matrix}\right)\neq (-1)^{k+1} e^{-\tpi \frac{a\lambda}{N}},
\label{eq:chi_minus}
\end{eqnarray}
we obtain
 \[
 \left[\left(D_z + \frac{2m\lambda}{\alpha} \right) \phi \left(z +\frac{\lambda}{\alpha}\tau +\frac{\mu}{\alpha},\tau \right)\right]_{z=0}
 =0.
 \]
Clearly $\chi\left(\begin{smallmatrix} -1&0\\0&-1\end{smallmatrix}\right)=\pm 1$.
If $\chi\left(\begin{smallmatrix} -1&0\\0&-1\end{smallmatrix}\right)=(-1)^{k+1}$ then \eqref{eq:chi_minus} fails for $\lambda =0$. Furthermore, if $\chi\left(\begin{smallmatrix} -1&0\\0&-1\end{smallmatrix}\right)=(-1)^k$ and $N$ is even then \eqref{eq:chi_minus} fails for $\lambda =N/2$. Hence the result holds provided $\chi\left(\begin{smallmatrix} -1&0\\0&-1\end{smallmatrix}\right)= (-1)^k$ and $N$ is odd,
completing the proof.
\hfill $\Qed$
\medskip

We find \eqref{eq:OnJacobiForms} does not necessarily hold for all $m$ and $\alpha$, however, as one would expect. For example, taking $m=1$ and $\lvert \alpha \rvert \geq 3$ shows that $\mathcal{T}_{k,\alpha}$ does not preserve $J_{k,1}$. This follows from the next lemma and the fact there are index $1$ Jacobi (cusp) forms such that $c(1,1) \not =0$ (see, for example, $\phi_{10,1}$ in \cite[display (17)]{EZ}). For the next result, we also note that $1+2\sqrt{m}(\sqrt{2}-1)>\sqrt{m}$ for all $m<17+12\sqrt {2}$ i.e., $m\le 33$.

\begin{lemma} \label{lem:JacobiCoeffs2}
Suppose $\phi \in J_{k,m}$ has  Fourier expansion \eqref{eq:phiFourier} with trivial multiplier and $\mathcal{T}_{k,\alpha} (\phi)$ has no poles. Set $h_m := \operatorname{max} \{t\in \mathbb{Z}: t\leq 2\sqrt{m} \} = \lfloor 2\sqrt{m} \rfloor$. If $\lvert \alpha \rvert \geq 1+2\sqrt{m}(\sqrt{2}-1)$ and $h_m \not = \frac{2m}{\lvert \alpha \rvert}$, then $c\left(1,h_m \right) =0$.
Furthermore, if $\lvert \alpha \rvert > \sqrt{m}$, then $c\left(1,h_m \right) =0$ in all cases.
\end{lemma}

\noindent \textbf{Proof.}
 Equation \eqref{eq:PoleProof1} of Lemma \ref{lem:LemPreserve} implies that in order for $\mathcal{T}$ to not introduce poles, we must have
 \[
 \sum_{\substack{n\geq 0 \\ r^2 \leq 4mn}} \left( r+\frac{2m\lambda}{\alpha} \right) c(n,r) e^{2\pi i\frac{\mu r}{\alpha}} q^{n+\frac{\lambda r}{\alpha}} =0,
 \]
 for all $\lambda, \mu \in \mathbb{Z}$
 Consider $n>0$ and pick $\lambda =\operatorname{sgn}(\alpha)$. The possible $q$ exponents are
 \[
 \left\{ n +r \frac{1}{\lvert \alpha \rvert} {\hspace{2mm} : \hspace{2mm}} 0\leq \lvert r\rvert  \leq \left\lfloor 2\sqrt{mn}\right\rfloor \right\}.
\]
This implies the smallest $q$ exponent for $n\ge 1$
is
\[
n - {\left\lfloor 2\sqrt{mn}\right\rfloor}\frac{1}{\lvert \alpha \rvert}\ge
n - 2\sqrt{mn}\frac{1}{\lvert \alpha \rvert}.
\]
For $n\ge 2$ we find
\[
\frac{2\sqrt{m}-1}{\sqrt{n}+1}+\frac{\sqrt{n}}{n-1}\le
\frac{2\sqrt{m}-1}{\sqrt{2}+1}+\sqrt{2}
=\left(2\sqrt{m}-1\right) \left(\sqrt{2}-1\right)+\sqrt{2}
\le \lvert \alpha \rvert,
\]
by assumption.
Hence
\[
(n-1)\lvert \alpha \rvert \geq 2\sqrt{m}\left(\sqrt{n}-1\right)+1,
\]
which implies that
\[
n - 2\sqrt{mn}\frac{1}{\lvert \alpha \rvert}\geq 1-\frac{\left(2\sqrt{m}-1\right)}{\lvert \alpha \rvert}.
\]
Finally, noting that $h_m>2\sqrt{m}-1$ we find
 \[
n - {\left\lfloor 2\sqrt{m
n}\right\rfloor}\frac{1}{\lvert \alpha \rvert} >  1-\frac{h_m}{\lvert \alpha \rvert},
 \]
i.e., the smallest $q$ exponent for $n\ge 2$ is bounded below by the smallest $q$ exponent for $n=1$.
The coefficient of this term is
 \begin{equation}
 \label{eq:h_mCoeff}
 \left( -h_m + \frac{2m}{\lvert \alpha \rvert} \right) c(1,-h_m).
 \end{equation}
 Since \eqref{eq:h_mCoeff} must equal $0$ if $\mathcal{T}_{k,m} (\phi)$ has no poles, we have $c(1,-h_m)=0$ unless $h_m = \frac{2m}{\lvert \alpha \rvert}$. Finally, we note $c(1,h_m)=(-1)^k c(1,-h_m)=0$ in this case.

If $\lvert \alpha \rvert>\sqrt{m}$ then since $h_m\le 2\sqrt{m}$ it follows that $h_m <\frac{2m}{\lvert \alpha \rvert}$ and  $c(1,h_m)=0$ in all cases.
  \hfill $\Box$\\

While the previous lemma is only relevant to Jacobi forms, there are many partition functions in the theory of VOAs known to be such (see examples in the next sections). The results in the next section, however, also pertain to the weak Jacobi forms arising in the theory of VOAs.

\section{Differential operators for strongly regular VOAs}\label{Sec:DiffOpsOnVOAs}

%
%

\subsection{Weak Jacobi forms and strongly regular VOAs}

A vertex operator algebra $V$ is said to be \textit{regular} if every weak $V$-module is a direct sum of irreducible ordinary modules. By \cite{ABD, Li-Finiteness}, this is equivalent to $V$ being rational and $C_2$-cofinite. This then implies that $V$ has finitely many inequivalent irreducible ordinary modules \cite{DLM-twisted,Zhu}, which we denote by $M^1 ,\dots ,M^r$. For an arbitrary such $V$-module, we often simply write $M$. A regular VOA which satisfies $L(1)V_1=0$ and is also of CFT-type, that is decomposes as $V=\mathbb{C} \vac\oplus \bigoplus_{n\geq 1} V_n$, is called \textit{strongly regular}.

Strongly regular VOAs come equipped with additional structure. For one, such VOAs have a nondegenerate symmetric invariant bilinear form $\langle \cdot ,\cdot \rangle \colon V\times V \to \mathbb{C}$ which is unique under the normalization $\langle \vac ,\vac \rangle =-1$ (see \cite{FHL,Li-Bilinear}). Additionally, the weight one space $V_1$ is a reductive Lie algebra and every $V$-module is completely reducible as a $V_1$-module \cite{DM-effective}. These facts are utilized below. For now, we focus on the fact that the space of $1$-point functions, $\operatorname{span}_\mathbb{C} \{V_{M^j}^J (v; z,\tau) \mid 1\leq j\leq r\}$, is closed under the the standard action of the Jacobi group. Indeed, let $J\in V_1$ be such that $J(0)$ acts semisimply on each $M^j$ with integral eigenvalues. Then the following theorem can be found in \cite{KrauelMasonI,KrauelMasonII}.

\begin{theorem}\label{thm:VOA-JacobiForms}
  Suppose $V$ is a strongly regular VOA. For each $v\in V_{[\wt [v]]}$ satisfying $J(n)v=0$ for $n\geq 0$ we have for all $\gamma = \left(\begin{smallmatrix} a&b \\c&d \end{smallmatrix}\right) \in \operatorname{SL}_{2}(\mathbb{Z})$ that
 \begin{equation}\notag %
 Z_{M^j}^J \left(v; \gamma .z,\gamma \tau \right) = (c\tau +d)^{\wt [v]} e^{\pi i \frac{\langle J,J\rangle c z^2}{c\tau +d}} \sum_{k =1}^r A_{j k}^\gamma Z_{M^k}^J (v;z,\tau),
 \end{equation}
 and for all and $[\lambda, \mu] \in \mathbb{Z}^2$ there exists a $j' \in \{1, \dots ,r\}$ such that
 \begin{equation}\notag %
 Z_{M^j}^J \left(v; z+\lambda \tau +\mu,\tau \right) = e^{-\pi i \langle J,J\rangle \left(\lambda^2 \tau +2\lambda z\right)} Z_{M^{j'}}^J (v;z,\tau).
 \end{equation}
\end{theorem}

Additionally, due to the $V$-module grading $M^j = \bigoplus_{n\geq 0} M^j_{n+\lambda_j}$, where $\lambda_j$ is the conformal weight of $M^j$, it immediately follows that each function $q^{\frac{\mathbf{c}}{24} - \lambda_j}Z_{M^j}^J \left(v;  z, \tau \right)$ has a Fourier expansion in a similar form to a weak Jacobi form. Indeed, vectors formed from the functions $Z_{M_j}^J$ ranging over the irreducible ordinary modules transform as vector-valued weak Jacobi forms. For an arbitrary element $v\in V$, however, it is unknown whether $Z_{M^j}^J \left(v; z, \tau \right)$ converges.


Before moving on, we discuss the relationship between the modes $J[n]$ and $J(n)$ given in \eqref{eq:asqround}.
\begin{lemma}
\label{lem:Jbrackets}
We have that $J(n)v=0$ for all $n\ge 0$ if and only if $J[0]v=J[1]v=0$.
\end{lemma}
\noindent \textbf{Proof.}
Recall (see, for example, \cite{Zhu}) that for a homogeneous element $v\in V$ we have \begin{equation}
 \notag 
 v[m] = m! \sum_{j\geq m} c(\wt (v),j,m)v(j).
\end{equation}
for  coefficients $c(\wt (v),j,m)$ defined via
\begin{equation}
 \notag 
 m!\sum_{j\geq m} c(\wt (v),j,m)x^j := \left(\text{Log} (1+x)\right)^m (1+x)^{\wt (v) -1}.
\end{equation}
Thus, in particular, we find since $\wt(J)=1$ that
\begin{align}
J[0]=J(0),\quad J[1] = \sum_{n\geq 1} \frac{(-1)^{n+1}}{n}J(n).
\label{eq:Jsq01}
\end{align}
 Therefore, $J[1]v=0$ if $J(n)v=0$ for all $n\geq 1$. Conversely, since $J(n) \colon M_{\lambda +k} \to M_{\lambda +k-n}$ then $J[1]v=0$ implies $J(n)v=0$ for all $n\geq 1$.     \hfill $\Box$\\


\subsection{Strongly regular VOAs with an $\slthat$ subalgebra}\label{Sec:Sl2}
\begin{proposition}\label{Prop:sltwo}
Suppose $V$ is a strongly regular VOA and that $\slt\subset V_1$ with standard generators $J,x,y$. Then the following are true:
\begin{enumerate}[label=\normalfont(\roman*)]
	\item $V$ contains an $\slthat$ Kac-Moody subalgebra of level $m=\half\langle J,J\rangle$;
	\item $J(0)$ has integer eigenvalues on all ordinary $V$-modules.
\end{enumerate}
\end{proposition}
\noindent{\textbf{Proof.}}
The standard basis $J,x,y$ satisfy the relations
\begin{align}
J(0)x=2x,\quad J(0)y=-2y,\quad x(0)y=J.
\label{eq:sl2}
\end{align}
Invariance of the bilinear form implies
$$
2\langle   x, y\rangle =\langle  J(0)x ,y \rangle
= -\langle  x(0) J ,y \rangle= \langle  J ,x(0) y \rangle =\langle   J ,J\rangle.
$$
Furthermore, $
\langle  x, J \rangle  = \langle  x ,x(0) y \rangle=-\langle  x(0)x , y \rangle=0$
and, similarly, $\langle  y, J \rangle=0$. In summary,
$J,x,y$ satisfy the identities
\begin{align*}
&\langle J ,J\rangle=2\langle x ,y\rangle,\quad
\langle x ,J\rangle=\langle y ,J\rangle=0.
\end{align*}
From the general VOA commutator formula we have
for all $u,v\in \slt$ that
\begin{align}
[u(r),u(s)]=(u(0) v)(r+s)+r\langle u ,v\rangle \delta_{r+s,0},
\label{eq:KM1}
\end{align}
where   $u(0)v=-v(0)u=[u,v]$ (by skew-symmetry).
For $ \langle J,J\rangle\neq 0$, $\langle \, ,\,\rangle$ is invertible and proportional to the standard $\slt$ invariant form $( \, ,\,)$ normalised  by $(J,J)=2$  so that
\begin{align*}
\label{invform}
\langle u,v\rangle=m(u,v),
\end{align*}
for $m=\half\langle J,J\rangle$,
which relation also holds for $m=0$.
Hence $J,x,y$  generate a $\slthat$ Kac-Moody algebra of level $m$ with relations \eqref{eq:KM1}.

Lastly, since $[u(0),L(0)]=0$ for all $u\in \slt$ then every $L(0)$ homogeneous space of a $V$-module  is a finite dimensional $\slt$-module. All finite dimensional $\slt$-modules are completely reducible into irreducible modules for which $J(0)$ has integral eigenvalues e.g.  Section~7 of \cite{Humph}. \hfill $\Qed$

\begin{corollary}\label{cor:slthat}
Suppose $V$ is strongly regular. If $V_1$ is non-Abelian then $V$ contains an $\slthat$ subalgebra of  positive integral level $m$.
\end{corollary}

 \noindent \textbf{Proof.} The space $V_1$ is a reductive Lie algebra by \cite{DM-effective}. Since $V_1$ is non-Abelian it must contain  an $\slt$ subalgebra. Furthermore, the $\slthat$ level $m$ is positive integral by Theorem~3.1 of \cite{DM-Integrability}.
\hfill $\Qed$

\medskip


\subsection{Distinguished degree $2$ differential operators}\label{subsect:diffops}
We now construct, in a VOA setting, differential operators of the type $\mathcal{M}$ defined in \eqref{eq:Moperator} in the case $\alpha= 2$.
By Propositions \ref{prop:apnpt} and \ref{prop:apnpt0}, for $v\in V$, we have
 \begin{equation} \label{eq:J[-1]Dz}
 Z_M^J (J[-1]v;z,\tau)  
 = D_z \left(Z_M^J (v;z,\tau)\right) +  \sum_{k\geq 2} G_{k} (\tau) Z_M^J (J[k-1]v;z,\tau),
 \end{equation}
using $D_z(\zeta^{J(0)})=J(0)\zeta^{J(0)}$.
 The operator $D_z$ is a well-known degree $1$ differential operator which preserves quasi-Jacobi forms (increasing their weight by $1$) but not (weak) Jacobi forms. We now turn to collecting degree two differential operators of (weak) Jacobi forms.

It is well-known that the $-1$ mode of  $\widetilde{\omega} \in V_{[2]}$ gives
 \begin{align*}
 Z_M^J (L[-2]v;z,\tau) &= \left(D_\tau +\wt[v]G_2 (\tau)\right) Z_M^J (v;z,\tau) +\sum_{k\geq 4} G_{k} (\tau) Z_M^J (L[k-2]v;z,\tau) \\
 &= {\vartheta}_{\wt[v]} \left(Z_M^J (v;z,\tau)\right) +\sum_{k\geq 4} G_{k} (\tau) Z_M^J (L[k-2]v;z,\tau).
 \end{align*}

 Multiples and powers of $J[-1]^2$ are considered in \cite{GK-differential} for $N=2$ superconformal field theories. Using Proposition \ref{prop:apnpt0} it is easily found that
\begin{align*}
 Z_M^J (J[-1]^2v;z,\tau) = &\left(D_z^2 +\langle J,J\rangle G_2 (\tau)\right) Z_M^J (v;z,\tau)
+2D_z\sum_{k\geq 2} G_{k} (\tau) Z_M^J (J[k-1]v;z,\tau) \\
 &  +\sum_{k,\ell \geq 2} G_{k}(\tau) G_{\ell} (\tau) Z_M^J (J[\ell -1]J[k-1]v;z,\tau).
\end{align*}
 The operator $\left(D_z^2 +\langle J,J\rangle G_2 (\tau)\right)$ also does not preserve Jacobi forms. However, a combination of it and ${\vartheta}$ does. Indeed, similarly to \cite[Subsection $3.1$]{GK-differential}, we find that
 \[
 Z_M^J \left(L[-2]\vac - \frac{1}{2\langle J,J\rangle} J[-1]^2 \vac;z,\tau \right) = \mathcal{H} \left(Z_M^J (z,\tau)\right),
 \]
  where $m=\langle J,J\rangle /2$ for the modified heat operator $\mathcal{H}$ of  \eqref{eq:operators}.
	Note that $J[0]v=J[1]v=0$ for $v=L[-2]\vac - \frac{1}{2\langle J,J\rangle} J[-1]^2 \vac$ in concurrence with Theorem~\ref{thm:VOA-JacobiForms}.

 We now consider other relevant elements which do not occur in \cite{GK-differential} but do occur for a strongly regular VOA $V$ with an $\slthat$ subalgebra.
In this case, we have $x,y\in V_{[1]}$ such that
\[
x [0] y   =J,\quad J[0]x=2 x,\quad J[0]y=-2 y,
\]
so that $J(0)$ has integral eigenvalues by Proposition~\ref{Prop:sltwo}.
Then for $v\in V_{[\wt [v]] }$ with $J[0]v=J[1]v=0$, Propositions \ref{prop:apnpt} and \ref{prop:apnpt0} give
\begin{equation}
\notag 
 \begin{aligned}
 Z_M^J (x[-1]y[-1]v ;z,\tau)
 =&  \left( - \widetilde{G}_{1}(2  z,\tau)D_z  + \langle x,y\rangle \widetilde{G}_{2}(2  z,\tau)  \right)Z_M^J (v;z,\tau) \\
 & +  \sum_{\ell\geq 2} (-1)^{\ell+1} \left(\widetilde{G}_{1}(2  z,\tau) {G}_{\ell}(\tau)+  \widetilde{G}_{\ell+1}(2  z,\tau) \right) Z_M^J (J[\ell-1]v;z,\tau) \\
 & +  \sum_{k ,\ell\geq 0} (-1)^{k +\ell} \widetilde{G}_{k +1}(-2  z,\tau)\widetilde{G}_{\ell+1}(2  z,\tau) Z_M^J ( y[k]x[\ell] v;z,\tau),
 \end{aligned}
\end{equation}
 where we use \eqref{eq:J[-1]Dz}. Noting that
\begin{equation*}
 \label{eq:component1}
 J[1]x[-1]y[-1]v =  -2  x[-1]y[0]v +2  y[-1]x[0] v +2  J[-1] v ,
\end{equation*}
 it is not surprising that the operator $ - \widetilde{G}_{1}(2  z,\tau)D_z + \frac{2m}{2 } \widetilde{G}_{2}(2  z,\tau)$ does not preserve the space of Jacobi forms.
However, using that
  \begin{align}
 J[1]L[-2]v = J[-1]v ,\quad
 J[1]J[-1]^2v =2\langle J,J \rangle J[-1]v,
\label{eq:J1Ln}
 \end{align}
we find, that if $Z_M^J(v,\tau, z)$ satisfies the Jacobi form functional equations, then for any $A ,B \in \mathbb{C}$ satisfying $A\not = 2B \langle J,J\rangle$, endomorphisms of the form
\begin{equation}
\label{eq:Melement}
 M= M_{(A ,B)} := L[-2] + \frac{1}{A-2B\langle J,J\rangle} \left(BJ[-1]^2 - \frac{A}{2 } x[-1]y[-1] \right),
\end{equation}
preserve these Jacobi form functional equations when applied to appropriate $v$ satisfying
\begin{align}
J[0]v= J[1]v= x[0]v=y[0]v=0,
\label{eq:Jxyv}
\end{align}
so that $J[1]Mv=0$.
 In other words, $M$ gives rise to the operator $\mathcal{M}= \mathcal{M}_{\wt [v],2}=\mathcal{M}_{(A,B),\wt[v], 2  ,m}$ (where $m=\langle J,J\rangle/2$ is the index) of \eqref{eq:Moperator} via
\begin{align}
 Z_M^J (Mv;z,\tau) &= \mathcal{M} \left(Z_M^J (v;z,\tau)\right) +\sum_{k\geq 2} G_{k} (\tau) Z_M^J (L[k-2]v;z,\tau)
\notag
 \\
 &\hspace{5mm} +\frac{B}{A-2B\langle J,J\rangle}\Big[\Big(2D_z \sum_{k\geq 4} G_{k} (\tau) Z_M^J (J[k-1]v;z,\tau)
\notag
 \\
 &\hspace{5mm} +\sum_{k,\ell \geq 4} G_{k}(\tau) G_{\ell} (\tau) Z_M^J (J[\ell -1]J[k-1]v;z,\tau)\Big) \Big]
\notag
 \\
 &\hspace{5mm} +  \frac{A}{2A-4B\langle J,J\rangle}\Biggl[\sum_{m\geq 2} (-1)^{m } \left(\widetilde{G}_{1}(2  z,\tau) {G}_{m}(\tau)+  \widetilde{G}_{m+1}(2  z,\tau) \right) Z_M^J (J[m-1]v;z,\tau)
\notag
 \\
 &\hspace{5mm} +  \sum_{k ,m\geq 0} (-1)^{k +m+1} \widetilde{G}_{k +1}(-2  z,\tau)\widetilde{G}_{m+1}(2  z,\tau) Z_M^J ( y[k]x[m] v;z,\tau)\Biggr].
\label{eq:Mv}
\end{align}
In particular, taking $v=\vac$ implies $ Z_M^J (M\vac;z,\tau) = \mathcal{M} \left(Z_M^J (z,\tau)\right)$.
In general,  \eqref{eq:Mv} preserves any Jacobi form transformation properties that might be satisfied by $Z_M^J (v;z,\tau)$ under appropriate conditions (such as in Theorem~\ref{thm:VOA-JacobiForms}). Thus we have the following.
\begin{proposition}
\label{prop:preserves}
 Suppose $V$ is a strongly regular VOA with an $\slthat$ subalgebra. Then for $v\in V_{[\wt [v]]}$ satisfying $J[0]v= J[1]v=x[0]v=y[0]v=0$, we find $\mathcal{M}_{\wt [v],2}$ preserves the Jacobi form functional equations of Theorem \ref{thm:VOA-JacobiForms} (adding $2$ to the weight) when applied to $Z_M^J (v;z,\tau)$.
\hfill $\Qed$
\end{proposition}


\subsection{Some applications}

Let $V_{n,\alpha}:=\{v\in V_n:J(0)v=\alpha v\}$ and denote its dimension by $\dim V_{n,\alpha}$.
\begin{proposition}
\label{prop:dimV}
 Suppose $V$ is a
strongly regular VOA with an $\slthat$ subalgebra. Then for any $\lambda \in\mathbb{Z}\setminus \{0\}$, we have
 \begin{equation}
 \label{eq:SumZero}
 \sum_{\substack{n\geq 1 \\ \lambda \mid 2 n}} \left( \frac{4 n}{\langle J,J\rangle\lambda^2} - 1\right)  \dim V_{n, \frac{2 n}{\lambda}} =
 \sum_{\substack{n\geq 1 \\ \lambda \mid 2 n}} \left(\frac{4 n}{\langle J,J\rangle\lambda^2} - 1 \right) (-1)^{ \frac{2n}{\lambda}} \dim V_{n, \frac{2 n}{\lambda}} =
1 .
 \end{equation}
\end{proposition}

\noindent \textbf{Proof.}
 By Propositions \ref{prop:apnpt} and \ref{prop:apnpt0} we know that $\mathcal{T}_{0,2} \left(Z_V^J (z,\tau)\right) = Z_V^J (M_{(1,0)}\vac;z,\tau)$ does not introduce poles. Therefore, we must have for $\lambda,\mu \in \Z$ that (recalling $m= \langle J,J\rangle /2$, cf. \eqref{eq:PoleProof1})
 \begin{align*}
 0=\left[\left(D_z - m\lambda \right) Z_V^J \left(z -\frac{\lambda}{2}\tau +\frac{\mu}{2},\tau \right)\right]_{z=0}
 &= \left[\left(D_z - m\lambda\right) q^{-\frac{\mathbf{c}}{24}}\sum_{\substack{n\geq 0 \\ r\in \mathbb{Z}}} \dim V_{n,r} e^{\pi i\mu r} q^{n- \frac{\lambda r}{2}} \zeta^{r} \right]_{z=0} \notag \\
 &= q^{-\frac{\mathbf{c}}{24}} \sum_{n\geq 0} \sum_{r\in \mathbb{Z}} \left(r-\frac{\langle J,J \rangle \lambda}{2}\right) \dim V_{n,r} e^{\pi i \mu r } q^{n- \frac{\lambda r}{2}}.
 \end{align*}
 Thus in particular the constant term (as $q$-expansions) of this expression must be zero, which occurs whenever $r=2 n/\lambda$ for $\lambda\neq 0$ and gives
 \[
 \sum_{\substack{n\geq 0 \\ \lambda \mid 2 n}} \left( \frac{2 n}{\lambda}-\frac{\langle J,J\rangle \lambda}{2}  \right) \dim V_{n, \frac{2 n}{\lambda}} e^{2\pi i \frac{n\mu }{\lambda}} =0.
 \]
 Recalling that $\dim V_{0,0}=1$ gives \eqref{eq:SumZero} taking $\mu =0$ or $1$.
\hfill $\Qed$
\medskip

We provide two simple corollaries that exploit the previous theorem.

\begin{corollary}
 Suppose $V$ satisfies the  conditions of Proposition~\ref{prop:dimV}. Then $V_{n,\pm 2 n} \not = \{0\}$ for finitely many $n\geq \lfloor \frac{\langle J,J\rangle}{4} +1 \rfloor$.
\end{corollary}

\noindent \textbf{Proof.}
 Taking $\lambda =\pm 1$  in \eqref{eq:SumZero} gives
 \[
 \sum_{n\geq 1} \left(\frac{4 n}{\langle J,J\rangle}-1 \right) \dim V_{n,\pm 2 n} = 1.
 \]
 We note that $4 n-\langle J,J\rangle>0$ if $n\geq \lfloor \frac{\langle J,J\rangle}{4} +1 \rfloor$,
$4 n -\langle J,J\rangle<0$ if $n\leq \lceil \frac{\langle J,J\rangle}{4} -1 \rceil$. Therefore, we can rewrite the sum above as
 \[
 \sum_{n\geq \left\lfloor \frac{\langle J,J\rangle}{4} +1 \right\rfloor} \left(\frac{4 n}{\langle J,J\rangle}-1 \right) \dim V_{n,  \pm 2 n} =
 1 + \sum_{n\geq 1}^{\left\lceil \frac{\langle J,J\rangle}{4} -1 \right\rceil} \left(1-\frac{4 n}{\langle J,J\rangle} \right) \dim V_{n,  \pm 2 n} \ge 0.
 \]
Thus $\dim V_{n,  \pm 2 n} \neq 0$ for only finitely many $n\geq \lfloor \frac{\langle J,J\rangle}{4} +1 \rfloor$.
\hfill $\Qed$
\medskip

The other corollary is the following.

\begin{corollary}
\label{cor:dimV1}
 Suppose $V$ satisfies the  conditions of Proposition~\ref{prop:dimV} and that $\langle J,J\rangle=2$. Then $\dim V_{1,\pm 2}=1$ and $V_{n,\pm 2 n}  = \{0\}$ for all $n\geq 2$.
\end{corollary}
\noindent \textbf{Proof.}
Taking $\lambda=\pm 1$ in \eqref{eq:SumZero} we find $\sum_{ n\geq 1 } \left(2 n - 1 \right)  \dim V_{n, \pm 2 n } =  1$ and the result follows. \hfill $\Qed$

%

\subsection{Strongly regular examples}

For more information and details on the following examples, we refer the reader to \cite{LL}, for example.

\subsubsection{Example 1: The VOA associated to the lattice $E_8$}

Throughout this section, let $V_{E_8}$ denote the VOA associated to the $E_8$ lattice $L_{E_8}$.
Then $V_{E_8}$ is a holomorphic VOA and
 \[
 V_{E_8}= \bigoplus_{\lambda \in L_{E_8}} M(1)\otimes e^{\lambda} .
 \]
 For any $a,b \in L_{E_8}$ recall that $a(0)b =0$, $a(0)e^b =\langle a,b \rangle e^b$, and $L(0)e^b =\frac{\langle b ,b \rangle}{2} e^b$.

 Assume that $h\in (V_{E_8})_1 =L_{E_8}$ has the property that $h(0)$ has integral eigenvalues on $V_{E_8}$ and $\langle h,h\rangle =2$, so that $h(1)h=\langle h,h\rangle\vac=2\vac$ and $m=\frac{\langle h,h\rangle}{2}=1$. We note that
 \begin{equation}
 \notag 
 Z_{V_{E_8}}^h(z,\tau)= \eta (\tau)^{-8}\sum_{\lambda \in L_{E_8}} \zeta^{\langle h,\lambda\rangle} q^{\frac{\langle \lambda ,\lambda\rangle}{2}} =\frac{\theta_{L_{E_8}}(z,\tau)}{\eta (\tau)^{8}}=\frac{E_{4,1}(z,\tau)}{\eta (\tau)^{8}},
 \end{equation}
 where $E_{4,1}$ is the weight $4$ index $1$ Jacobi-Eisenstein series (see \cite[Section $2$, equation ($1$)]{EZ}). Thus, $\eta (\tau)^8 Z_{V_{E_8}}(z,\tau)$ is a Jacobi form of weight $4$ and index $\frac{\langle h,h\rangle }{2}=1$.

Consider also the elements $e^h ,e^{-h} \in (V_{E_8})_1$. Note that for $n\geq 0$
$$
h(n)e^{\pm h} =h[n]e^{\pm h} =\delta_{n,0} \langle h,\pm h \rangle e^{\pm h} =\pm 2e^{\pm h}
$$
and
  \[
  e^h [n] e^{-h} =\begin{cases} h &\operatorname{ if } n=0, \\ \vac &\operatorname{ if } n=1,\\ 0&\operatorname{ otherwise.} \end{cases}
  \]
  Recalling Propositions \ref{prop:apnpt} and \ref{prop:apnpt0} with $x=e^h$ and $y=e^{-h}$, we find that
  \begin{equation}
  \notag %
  Z_{V_{E_8}}^h \left( e^h [-1] e^{-h};z,\tau \right) =\left(-\widetilde{G}_{1}(2 z,\tau)D_z +\widetilde{G}_{2}(2 z,\tau) \right)Z_{V_{E_8}}^Jh (z,\tau),
  \end{equation}
  so that
  \begin{align*}
  Z_{V_{E_8}}^h  \left( T\vac ;z,\tau \right) &=\mathcal{T}_{0,2} \left(Z_{V_{E_8}}^h  \left(z,\tau \right) \right)= \mathcal{T}_{0,2}  \left( \frac{E_{4,1}(z,\tau)}{\eta (\tau)^{8}} \right), \\
  Z_{V_{E_8}}^h  \left( S\vac ;z,\tau \right) &=\mathcal{S}_{0,2} \left(Z_{V_{E_8}}^h
	\left(z,\tau \right)\right) = \mathcal{S}_{0,2}  \left( \frac{E_{4,1}(z,\tau)}{\eta (\tau)^{8}} \right),
  \end{align*}
	where $S=M_{(0,1)}$ and $T=M_{(1,0)}$. 	
   Moreover, since $h[1]S\vac=h[1]T\vac=0$, we have $\mathcal{S}_{0,2} ( \frac{E_{4,1}(z,\tau)}{\eta (\tau)^{8}} )$ and $\mathcal{T}_{0,2} ( \frac{E_{4,1}(z,\tau)}{\eta (\tau)^{8}} )$ transform like Jacobi forms of weight $2$ and index $1$.

 Indeed, it can be found
$$
\mathcal{S}_{0,2} \left(Z_{V_{E_8}}^h (z,\tau)\right) = -\frac{7}{24} \frac{E_{6,1} (z,\tau)}{\eta (\tau)^{8}}, \qquad
  \mathcal{T}_{0,2} \left(Z_{V_{E_8}}^h (z,\tau)\right) = -\frac{7}{24} \frac{E_{6,1} (z,\tau)}{\eta (\tau)^{8}}.
	$$
 Such expressions are similar to one of the three Ramanujan equations studied in \cite[Corollary $3$]{Oberdieck-Serre}.


\subsubsection{Example 2: The VOA associated to affine Lie algebra $\widehat{\frak{sl}}_2$}

\indent Consider the (strongly regular) simple VOA $V:= L_{\widehat{\frak{sl}}_2}(m ,0)$ associated to the affine Lie algebra $\widehat{\mathfrak{sl}}_2$ of level $m \in \mathbb{N}$, where $h,x,y \in \mathfrak{sl}_2 =V_1$ are the typical basis elements of the Lie algebra. Such VOAs have $m+1$ many inequivalent irreducible modules, which we denote here as $V=M^0 ,M^1, \dots ,M^m$. Then we have (since $h(0)x=[h,x]=2x$ and $\langle x,y\rangle =m$) that for every $M^j$
 \begin{equation}
 \notag %
 Z_{M^j}^h \left( T\vac;z,\tau \right) =\mathcal{T}_{0,2} \left( Z_{M^j}^h (z,\tau) \right) \quad \text{and} \quad
 Z_{M^j}^h \left( S\vac;z,\tau \right) =\mathcal{S}_{0,2} \left( Z_{M^j}^h (z,\tau) \right) ,
 \end{equation}
	where as before $S=M_{(0,1)}$ and $T=M_{(1,0)}$.
 Since $h[1]S\vac=h[1]T\vac=0$, we again have that
  \[
 \left(Z_{M^0}^h \left(T\vac;z,\tau\right) ,\dots ,Z_{M^m}^h \left(T\vac;z,\tau\right)\right)^{T} \quad \text{and} \quad
 \left(Z_{M^0}^h\left(S\vac;z,\tau\right) ,\dots ,Z_{M^m}^h \left(S\vac;z,\tau\right)\right)^{T}
 \]
transform as vector-valued weak Jacobi forms of weight $2$ and index $\frac{\langle h,h\rangle}{2}=m$.

 Taking $m=1$ (that is, the case of $A^{(1)}_1$), it is known that $\mathbf{c}=1$ and
 \[
 Z_V^h (z,\tau) = \frac{\sum_{n\in \mathbb{Z}} \zeta^{n}q^{n^2}}{\eta(\tau)} =: \frac{\theta_3 (z,2\tau)}{\eta (\tau)},
\qquad
Z_{M^1}^h (z,\tau) = \frac{\sum_{n\in \mathbb{Z}} \zeta^{n+\frac{1}{2}}q^{\left(n+\frac{1}{2}\right)^2}}{\eta(\tau)} =: \frac{\theta_2 (z,2\tau)}{\eta (\tau)},
 \]
 (for example, these can be deduced from \cite{Kac-infinite}).
Therefore,
 \[
 \begin{pmatrix}
 Z_{V}^h (S\vac;z,\tau) \\ Z_{M^1}^h  (S\vac;z,\tau)
 \end{pmatrix}
 =
\begin{pmatrix}
  \mathcal{S}_{0,2} \left(Z_{V}^h  (z,\tau)\right) \\ \mathcal{S}_{0,2} \left(Z_{M^1}^h  (z,\tau)\right)
 \end{pmatrix}
 =
 \frac1{\eta (\tau)} \begin{pmatrix}
 \mathcal{S}_{\frac{1}{2},2} \left(\theta_3 (z,2\tau)\right) \\ \mathcal{S}_{\frac{1}{2},2} \left(\theta_2 (z,2\tau)\right)
 \end{pmatrix},
 \]
 where we extended the definition of $\mathcal{S}_{k,m}$ to $k\in \mathbb{Q}$ in the natural way.

 Finally, we note that the important result that $x[-1]^k \vac =0$ for $k\geq 2$ is known due to nilpotency arguments (see, for example, \cite{LL}). However, this is also now immediate from Corollary \ref{cor:dimV1}.


\section{Fermionic models}\label{Sec:Fermion}
\subsection{Vertex operator super algebras}
In this section we consider an analogue of the $\slt$ structure of Section~\ref{Sec:DiffOpsOnVOAs} for a central charge $\mathbf{c}$ vertex operator super algebra (VOSA) of  CFT-type $V=\oplus_{k \in \half\Z}V_{k}$ e.g. \cite{MTZ}.  We define a parity operator $p(v)\equiv  2k\pmod{2}$ for $v\in V_{k}$ and define a ``fermion number'' automorphism $\sigma$ by $\sigma v :=(-1)^{p(v)}v$.

Assume that there exists $2R$ ``free fermion'' vectors  $\psi_r^{\pm}\in V_{\half}$ for $r=1,\ldots ,R$ with vertex operators
$Y\left(\psi_r^{\pm },z \right)=\sum_{n\in \mathbb{Z}}\psi_r^{\pm }(n)z^{-n-1}$ such that $\psi_r^{+}(0)\psi_s^{-}=\delta_{rs}\vac$ and $\psi_r^{\pm}(0)\psi_s^{\pm}=0$. This implies the anti-commutator relations
\begin{align}
\left[\psi_r^{+}(m),\psi_s^{-}(n)\right]=\delta_{rs}\delta _{m,-n-1},\quad
\left[ \psi_r^{\pm}(m),\psi_s^{\pm }(n)\right]=&0. \notag 
\end{align}
Defining
\begin{align*}
J:=\sum_{r=1}^R \psi_r^{+}(-1)\psi_r^{-},
\end{align*}
it follows  (analogously to \eqref{eq:sl2}) that
 \begin{align*}
J(0)\psi_r^{\pm}=\pm \psi_r^{\pm},\quad J(1)J =R\vac,
\end{align*}
i.e.,  $\langle J,J\rangle =R$ with
\begin{align}
\left[J(m),\psi_r^{\pm}(n)\right]=\pm\psi_r^{\pm}(m+n),\quad [ J(m),J(n)]=mR\delta _{m,-n}.
\label{eq:Jpsicom}
\end{align}
We further assume that the fermion number automorphism is given by
\begin{align}
\sigma=e^{\pi i J(0)},
\label{eq:sigmaJ}
\end{align}
so that $J(0)$ has integral eigenvalues on $V$.
\medskip

In order to illustrate the main results of this paper, we consider the $\sigma$-twisted $V$-module (the Ramond sector).
Define for all $v\in V$
\[
Y_{\sigma}(v,z):=Y(\Delta (\sigma, z)v,z), 
\quad \Delta(\sigma,z):=z^{\half J(0)}\exp \left(- \half\sum\limits_{n\geq 1}\frac{J(n)}{n}%
(-z)^{-n}\right) .
\]
Then $(V,Y_{\sigma})$ is the $\sigma$-twisted $V$-module by a theorem of Li \cite{Li}.
In particular
\begin{align*}
Y_{\sigma}\left(\psi^\pm,z\right) &=z^{\pm \half}Y(\psi^\pm,z),\quad
Y_{\sigma}(J,z) =Y(J,z) +\half z^{-1}\id_V,
\\
Y_{\sigma}(\omega,z) &=Y(\omega,z) +\half z^{-1}Y(J,z) +\frac{R
}{8}z^{-2}\id_V,
\end{align*}
where $\omega$ is the Virasoro vector.
With $Y_{\sigma}(v,z)=\sum_{n\in \Z+\half p(v)}v_\sigma(n)z^{-n-1}$, we thus find that
\begin{align}
J_{\sigma}(0)=J(0)+\half,
\notag 
\quad
L_{\sigma}(0)=L(0)+\half J(0)+\frac{R}{8}.
\notag 
\end{align}

\subsection{Jacobi $n$-point functions}
We  define a Jacobi $0$-point function with half integral grading given by
\begin{align}
Z_V^J(z,\tau):=\str_{V}\zeta^{J(0)}q^{L(0)-\frac{\mathbf{c}}{24}},
\label{eq:ZVfermion}
\end{align}
for supertrace $\str_V A q^{L(0)}:=\tr_{V} \sigma A q^{L(0)}$ for $A\in\End V$.
The supertrace \eqref{eq:ZVfermion} associated with the Ramond $\sigma$-twisted module is
\begin{align} \label{eq:Zsig}
Z^J_{V_\sigma}(z,\tau):=
\tr_{V} e^{i\pi J_\sigma(0)} \zeta^{J_\sigma(0)}q^{L_\sigma(0)-\frac{\mathbf{c}}{24}}
=i\str_{V} \zeta^{J(0)+\half}q^{L(0)+\half J(0)-\frac{(\mathbf{c}-3R)}{24}}.
\end{align}
We can generalize \eqref{eq:Zsig} to all $\sigma$-twisted Jacobi $n$-point functions such as in \eqref{eq:npt}. These can be computed in terms of appropriate $n$-point functions on $V$ with a shifted  Virasoro vector \cite{MTZ}
\begin{align}
\omega_s=\omega-\half J(-2)\vac,
\label{eq:oms}
\end{align}
for  central charge $\mathbf{c}_s=\mathbf{c}-3R$ and  shifted grading operator
\begin{align*}
L_s(0)=L(0)+\half J(0).
\end{align*}
We note that $L_s(0)$ necessarily has integral eigenvalues on $V$ from \eqref{eq:sigmaJ}.
As shown in Proposition~9 of \cite{MTZ}, every  $\sigma$-twisted $n$-point function
\begin{align}
&Z_{V_\sigma}^J\left ((v_1,x_1),\dots, (v_n,x_n);z,\tau\right)
\notag
\\
& :=i\str_V
Y_\sigma(e^{x_1L(0)}v_1,e^{x_1}) \cdots Y_\sigma(e^{x_nL(0)}{v}_n,e^{x_n})
\zeta^{J(0)+\half}q^{L(0)+\half J(0)-\frac{(\mathbf{c}-3R)}{24}},
\label{eq:Zsignpt}
\end{align}
 can be expressed in terms of an appropriate untwisted $n$-point function with shifted Virasoro vector
\begin{align}
& Z_{V_\sigma}^J\left ((v_1,x_1),\dots, (v_n,x_n);z,\tau\right)
=
i \str_V Y(e^{x_1 L_s(0)}U{v}_1,e^{x_1}) \cdots
Y(e^{x_n L_s(0)}U{v}_n,e^{x_n})
\zeta^{J(0)+\half}q^{L_s(0)-\frac{\mathbf{c}_s}{24}},
\label{eq:Znpts}
\end{align}
 where (using \eqref{eq:Jsq01})
\begin{align*}
U:=\exp\left(\half \sum_{n\ge 1}\frac{(-1)^{n+1}}{n}J(n)\right)=e^{\half J[1]}.
\end{align*}
Since the shifted grading is integral we may apply the Zhu reduction formulas of Propositions~\ref{prop:Zhured}--\ref{prop:apnpt0} to the supertrace \eqref{eq:Znpts} taking due regard to the anti-commuting properties of fermion vertex operators\footnote{
The supertrace is required in \eqref{eq:Zsignpt} and \eqref{eq:Znpts} in order to obtain the appropriate quasi-Jacobi coefficient functions appearing in  Propositions~\ref{prop:Zhured}--\ref{prop:apnpt0}.}.
\subsection{A degree 1 differential operator}
Consider $v\in V_{[k]}$ such that
\begin{align}
J[0]v= J[1]v= \psi_r^{\pm}[0]v =0,
\label{eq:Jpsiv}
\end{align}
for $r=1,\ldots,R$
in analogy to \eqref{eq:Jxyv}. Using \eqref{eq:Jpsicom} we find that for the VOSA endomorphism\footnote{The superscript $S$ indicates  a super algebra endomorphism.}
\[
N^{S}:=J[-1]-\sum_{r=1}^R \psi_r^{+}[-1]\psi_r^{-}[-1],
\]
 then
$
J[1]N^{S}v=0$.
Furthermore, Zhu reduction using \eqref{eq:Znpts} implies
\begin{align*}
 Z_{V_\sigma}^J \left(N^{S}v;z,\tau\right)
 = &
\left( D_z +R\widetilde{G}_{1}( z,\tau)\right)Z_{V_\sigma}^J (v;z,\tau)
+  \sum_{k\geq 1} G_{2k} (\tau) Z_{V_\sigma}^J (J[2k-1]v;z,\tau)\\
 &+\sum_{k ,\ell\geq 1} (-1)^{k +\ell} \widetilde{G}_{k +1}(  z,\tau)\widetilde{G}_{\ell+1}( -z,\tau)
Z_{V_\sigma}^J \left( \sum_{r=1}^R\psi_r^{+}[k]\psi_r^{-}[\ell]v;z,\tau\right).
\end{align*}
In particular, for $v=\vac$ we find that $N^{S}\vac=0$ so that
\begin{align*}
\left( D_z + R\widetilde{G}_{1}( z,\tau)\right)Z_{V_\sigma}^J ( z,\tau) =0.
\label{eq:DzZ}
\end{align*}
Hence,  since $\widetilde{G}_{1}( z,\tau)=-P_1(z,\tau)=-D_z(\text{Log} (K(z,\tau)))$, we have
\begin{align*}
Z_{V_\sigma}^J ( z,\tau) =\theta_{1}(z,\tau)^R F(\tau),
\end{align*}
for Jacobi theta function \eqref{eq:theta1} and  some $z$-independent function $F(\tau)$.
This condition severely restricts the possible VOSAs satisfying \eqref{eq:sigmaJ}.  One obvious family of examples is the VOSA formed by taking the tensor product of the VOSA generated by $\psi^{\pm}_r$ and an arbitrary VOA $W$ for which
\[
Z_{V_\sigma}^J ( z,\tau) =\left(\frac{\theta_{1}(z,\tau)}{\eta(\tau)}\right)^R Z_W(\tau).
\]

\subsection{Degree 2 differential operators}
Much as in Section~\ref{subsect:diffops}, we consider differential operators of the type $\mathcal{M}$ defined in \eqref{eq:Moperator} for $\alpha= 1$ that arise in the Zhu reduction of appropriate VOSA $1$-point functions. These examples differ from those considered in \cite{GK-differential} for $N=2$ superconformal algebras. Consider $v\in V_{[k]}$ obeying \eqref{eq:Jpsiv}. Note that in addition to  \eqref{eq:J1Ln}  we have, using \eqref{eq:Jpsicom}, that
\begin{align*}
J[1] \psi^{+}_r[-2]\psi^{-}_r[-1]v&=\psi^{+}_r[-1]\psi^{-}_r[-1]v,
\\
J[1] \psi^{+}_r[-1]\psi^{-}_r[-1]J[-1]v&=\left(J[-1]+(R-2)\psi^{+}_r[-1]\psi^{-}_r[-1]\right)v.
\end{align*}
Thus for the VOSA endomorphism (cf. \eqref{eq:Melement})
\begin{multline*}
M^{S}=M^{S}_{(A,B)}\\
:=L[-2]
+ \frac{1}{A-2BR} \Big (BJ[-1]^2  + \frac{A}{R } \sum_{r=1}^R
\left((R-2)\psi_r^{+}[-2]\psi_r^{-}[-1]
-\psi_r^{+}[-1]\psi_r^{-}[-1]J[-1] \right)\Big),
\end{multline*}
we find
$
J[1]M^{S}v=0$. Zhu reduction  implies  the operator $\mathcal{M}=\mathcal{M}_{(A,B),\wt[v], 1 ,\half R}$ of  \eqref{eq:Moperator} with $\alpha=1$ occurs as follows:
\begin{align}
 Z \left(M^{S}v\right)
 = &
\mathcal{M}\left(Z (v)\right)  +  \frac{A}{A-2BR}
\left(2\widetilde{G}_{2}( z,\tau)+\widetilde{G}_{1}( z,\tau)^2\right)Z (v)
+\sum_{k\geq 4} G_{k} (\tau) Z (L[k-2]v)
\notag
\\
 &  +\frac{B}{A-2BR}\left[2D_z \sum_{k\geq 4} G_{k} (\tau) Z (J[k-1]v)
  +\sum_{k,\ell \geq 4} G_{k}(\tau) G_{\ell} (\tau) Z (J[\ell -1]J[k-1]v)\right]
\notag
\\
 &  +  \frac{A}{A-2BR}\left[
\sum_{m\geq 4}  \widetilde{G}_{1}( z,\tau) {G}_{m}(\tau) Z \left( J[m-1]v\right)\right.
\notag
\\
& +\frac{R-2}{R}  \sum_{k ,\ell\geq 1} \sum_{r=1}^R(-1)^{k +\ell} (k+1)\widetilde{G}_{k +2}( z,\tau)\widetilde{G}_{\ell+1}(  -z,\tau)  Z  \left( \psi^{+}_r[k]\psi^{-}_r[\ell] v\right)
\notag
 \\
 &
-\frac{1}{R}  \sum_{k ,\ell\geq 1} \sum_{r=1}^R(-1)^{k +\ell}  \widetilde{G}_{k +1}( z,\tau)\widetilde{G}_{\ell+1}(  -z,\tau)
\Big\{
\left(D_z + 2\widetilde{G}_{1}( z,\tau)\right)Z  \left( \psi^{+}_r[k]\psi^{-}_r[\ell] v\right)
\notag
\\
 & \qquad \qquad\qquad\qquad+Z  \left( \psi^{+}_r[k]\psi^{-}_r[\ell-1] v\right)
- Z  \left( \psi^{+}_r[k-1]\psi^{-}_r[\ell] v\right)
\Big\}
\notag
\\
&\left. +
 \frac{1}{R}  \sum_{k ,\ell\geq 1 }\sum_{r=1}^R \sum_{ m\ge 4}(-1)^{k +\ell}  \widetilde{G}_{k +1}( z,\tau)\widetilde{G}_{\ell+1}(  -z,\tau)  G_m(\tau)
Z \left( J[m-1]\psi^{+}_r[k]\psi^{-}_r[\ell] v\right)
\right],
\label{eq:Msv}
\end{align}
where, for notational simplicity we  abbreviate
$Z(u):=Z_{V_\sigma}^J (u;z,\tau)$  for $u\in V$. Equation \eqref{eq:Msv} is an $\alpha=1$ analogue of \eqref{eq:Mv}.
In particular, taking $v=\vac$ and using \eqref{eq:G_2tilde} we have
\begin{align}
Z_{V_\sigma}^J\left(M^{S}\vac;z,\tau\right)
 &=
\left [\mathcal{M} -  \frac{A}{A-2BR}
\frac{1}{4 \pi^2}\wp(z,\tau)\right]Z_{V_\sigma}^J (z,\tau).
\label{eq:MsZpart}
\end{align}

\subsection{The free fermion model }
We now specialize to the rank $2R$ free fermion VOSA $V$ generated by $\psi^{\pm}_r$ for $r=1,\ldots, R$ for which the $\sigma$-twisted supertrace is (e.g.\ \cite{MTZ})
\begin{align*}
Z_{V_\sigma}^J(z,\tau)=\left(\frac{\theta_1(z,\tau)}{\eta(\tau)}\right)^R,
\label{eq:Zsigfermion}
\end{align*}
for Jacobi theta function \eqref{eq:theta1}.

Consider first of all the case  $R=1$ and let $\psi^{\pm}=\psi_{1}^{\pm}$. Noting that
 $J(0)\psi^{+}=\psi^{+}$, then Corollary~\ref{cor:ZeroRes} implies that for all $z=\lambda \tau +\mu \in {\mathbb{Z}\tau} +\mathbb{Z}$, we find that
\begin{align*}
0=\sum_{m\ge 0}
\frac{\lambda^m}{m!}
Z_{V_\sigma}^J\left (\psi^{+}[m]\psi^{-};z,\tau\right)=
Z_{V_\sigma}^J\left (z,\tau\right),
\end{align*}
i.e., $\theta_1(\lambda \tau +\mu,\tau)=0$ for all $\lambda ,\mu\in\Z$, as is well-known.
The central charge $-2$ shifted Virasoro vector \eqref{eq:oms} in the square bracket formalism  is given by
\[
\widetilde{\omega}=L[-2]\vac=\half (J[-1]^2-J[-2])\vac=\psi^{+}[-2]\psi^{-}.
\]
Thus $M^S\vac=\frac{B}{A-2B}L[-1]J$ implying $Z_{V_\sigma}^J\left (M^S\vac ;z,\tau\right)=0$.
Hence, using \eqref{eq:Heat}, \eqref{eq:Top} and $D_\tau\left( \eta(\tau)\right)=-\half G_2(\tau)\eta(\tau)$ we find \eqref{eq:MsZpart} is equivalent to
\begin{align}
\label{eq:heat}
\mathcal{H}_{\half}\left(\theta_1(z,\tau)\right) &= \left(D_\tau -\half D_z^2\right)\theta_1(z,\tau)=0,
\\
\mathcal{T}_{\half,1}\left(\theta_1(z,\tau)\right) &= \frac{1}{4 \pi^2}\wp(z,\tau)\theta_1(z,\tau),
\label{eq:TZ}
\end{align}
for index $\half$.
 Display \eqref{eq:heat} is the classic heat equation whereas,
using $\widetilde{G}_{1}( z,\tau)=-P_1(z,\tau)$ and   $P_2(z,\tau)=-D_z \left(P_1(z,\tau)\right)$, we find that \eqref{eq:TZ} implies \eqref{eq:G_2tilde}.

We next  consider the $\sigma$-twisted $2n$-point function for $n$ copies of $\psi^{+}$ and $\psi^{-}$.
The first Zhu reduction formula Proposition~\ref{prop:Zhured} implies (cf.\ Proposition~14 of \cite{MTZ}) the next result.

\begin{proposition} We have
\label{prop:ferm}
\begin{align}
&Z_{V_\sigma}^J\left ((\psi^{+},x_1),(\psi^{-},y_1),\dots, (\psi^{+},x_n),(\psi^{-},y_n);z,\tau\right)
=
\det \left(\mathbf{P}\right)  \frac{\theta_1(z,\tau)}{\eta(\tau)},
\label{eq:nptdet}
\end{align}
for a matrix  $\mathbf{P}$ with components
$\mathbf{P}_{jk}=\widetilde{P}_{1}(x_j-y_k, z,\tau)$
for $1\le j,k\le n$.
\hfill $\Qed$
\end{proposition}

As demanded by Proposition~\ref{prop:Zhured0}, we find the following corollary.

\begin{corollary}
The $2n$-point function  \eqref{eq:nptdet} is convergent  for all $x_j\neq y_k$ with $x_j,y_k,z\in\C$ and $0<|q|<1$ .
\label{cor:nptdet}
\end{corollary}

\noindent \textbf{Proof.}
From Proposition~\ref{prop:P1tilde}~(v), we have simple pole structure
\begin{align*}
\mathbf{P}_{jk}=\frac{q_{y_k}^\lambda q_{x_j}^{-\lambda}}{\tpi(z-\lambda \tau -\mu)}+O(1),
\end{align*}
in the neighborhood of $z=\lambda \tau+\mu$
for each $\lambda ,\mu \in \Z$.
The singular parts of the columns of $\mathbf{P}$ are linearly dependent and hence $\det (\mathbf{P})$ has a pole of order $1$ for all $z\in {\mathbb{Z}\tau} +\Z$. Since $\theta_1(z,\tau)$ has a simple zero at $z\in {\mathbb{Z}\tau} +\Z$ the result follows.
\hfill $\Qed$
\medskip

Equation \eqref{eq:nptdet} is a generating function for all $\sigma$-twisted $n$-point functions as explained in Proposition~15 of \cite{MTZ}. Thus we conclude the convergence of $n$-point functions.

\begin{proposition} We have that
\label{prop:nptcgt}
$Z_{V_\sigma}^J\left ((v_1,x_1),\dots, (v_n,x_n);z,\tau\right)$ is convergent  for all $x_j\neq x_k$ with $j\neq k$ for  $x_j,z\in\C$ and $0<|q|<1$.
\end{proposition}

\medskip

We next consider the rank $2R$ fermion VOSA  generated by $\psi^{\pm}_r$ for $r=1,\dots ,R$  with shifted Virasoro vector $\widetilde{\omega}=\sum_{r=1}^R\psi^{+}_{r}[-2]\psi^{-}_{r}$ of central charge $-2R$. We can construct all
weak Jacobi forms (with the same multiplier system as $\theta_1(z,\tau)^R$) of the form $F_R(z,\tau)\theta_1(z,\tau)^R$ for all $R\ge 2$ for the meromorphic Jacobi form $$
F_R(z,\tau):=\frac{(-1)^{R+1}}{R}\left(P_R(z,\tau)-G_R(\tau)\right)
$$
of \eqref{eq:gen1} or $F_R(z,\tau)=K_R(z,\tau)$ of \eqref{eq:Kn} in terms of Jacobi 1-point functions for specific vectors in the kernel of $J[1]$.

Define for all $k\in \mathbb{N}_0$ and $1\le r\le  R$ the following commuting operators
\begin{align}
b_r^{k}:=
\begin{cases}
1&{\rm if }\ k=0,\\
\psi_r^{+}[-k]\psi^{-}_r[-1] &{\rm if }\ k\ge 1.
\end{cases}
\notag 
\end{align}
Noting that $(b^k_r)^2=0$ for all $k>0$, we consider the Fock vector
\begin{align}
b_{r_1}^{k_1}\cdots b_{r_m}^{k_m}\vac,
\label{eq:Bvect}
\end{align}
for $k_1,\dots,k_m\in \mathbb{N}$ with non-repeating fermion labels $1 \le r_1<r_2< \dots <r_m\le R$, i.e., any $k_j$ value occurs at most $R$ times.
The Fock vector \eqref{eq:Bvect} is of $L[0]$ weight $n=\sum_{j=1}^m k_j$, a partition of $n$ with $m$ parts.
It is useful to denote this partition by $\pi=1^{j_1} \dots n^{j_n} $
indicating that there are  $m=j_1+\dots +j_n$ parts of $n= j_1+2j_2+\dots +n j_n$ (but with $k$ occurring $j_k\le R$ times).

Define an  $L[0]$ weight $n$ vector for a partition $\pi=1^{j_1} \dots n^{j_n}$ of $n\le R$ by
\begin{align}
\beta^n_\pi:=\sum_{1 \le r_1< \dots <r_m\le R} b_{r_1}^{k_1}\cdots b_{r_m}^{k_m}\vac,
\label{eq:Bbold}
\end{align}
where the sum is taken over $\binom{R}{j_0j_1\dots j_R}$  independent vectors of the form \eqref{eq:Bvect},  where we  define $ j_0=R-m\ge 0$. For example, $\beta^3_{1^2}=(b_1^1b_2^1+b_1^1b_3^1+b_2^1b_3^1)\vac$.
\begin{proposition}
\label{prop:J1Bbold}
For $\pi=1^{j_1} \dots R^{j_R}$, a partition of $R$, we have
\begin{align*}
J[1]\beta^R_\pi = \sum_{k=1}^R
(j_{k-1}+1)\chi_{j_k}\beta^{R-1}_{\pi_k},
\end{align*}
where
$\pi_k = 1^{j_1}\dots (k-1)^{j_{k-1}+1}k^{j_k-1}\dots R^{j_R}$ is a partition of $R-1$ (provided $j_k>0$) and
where $\chi_0=0$ and $\chi_j=1$ for $j>0$.
\end{proposition}

\noindent \textbf{Proof.}
 From \eqref{eq:Jpsicom} we find that
\begin{align*}
J[1] b_{r_1}^{k_1}\dots b_{r_m}^{k_m}\vac =
\sum_{j=1}^{m}b_{r_1}^{k_1}\cdots b_{r_j}^{k_j-1}\cdots b_{r_m}^{k_m}\vac.
\end{align*}
Thus, provided that $j_k>0$, every independent summand vector in $\beta^{R-1}_{\pi_k}$ (cf.\ \eqref{eq:Bbold}) arises in $J[1]\beta^R_{\pi}$  with multiplicity $j_{k-1}+1$ which is the number of ways that a given set of $j_{k-1}+1$ fermion labels can be constructed from a subset of $j_{k-1}$  labels together with the remaining label.
\hfill $\Qed$
\medskip

We now define a vector of $L[0]$-weight $R$ in the kernel of $J[1]$, namely
\begin{align*}
\Phi^R:=\sum_{\pi}\frac{(-1)^{m+1}}{m\binom{R}{m}}\beta^R_{\pi},
\end{align*}
where the sum runs over all partitions $\pi=1^{j_1}\dots R^{j_R}$ of $R$ with $m=\sum_{k=1}^{R}j_k$ and $j_0=R-m$. For example, for $R=1,2$ we have
\begin{align*}
\Phi^1&=J,\quad \Phi^2=\frac{1}{2}\left(\psi_1^{+}[-1]\psi_1^{-}[-1]\psi_2^{+}[-1]\psi_2^{-}[-1] -\psi_1^{-}[-2]\psi_1^{+}[-1]-\psi_2^{-}[-2]\psi_2^{+}[-1]\right)\vac.
\end{align*}
We find the following result.
\begin{proposition} We have that
$J[1]\Phi^R=0$ for all $R\ge 2$.
\label{prop:J1ker1}
\end{proposition}

\noindent \textbf{Proof.}
Proposition~\ref{prop:J1Bbold} implies that $J[1]\Phi^R$ is a linear combination of independent vectors of the form $\beta^{R-1}_{\kappa}$ for a partition $\kappa=1^{j_1}\dots (R-1)^{j_{R-1}}$
 of $R-1$ into $m=j_1+\cdots +j_{R-1}$ parts. Each such vector appears in $J[1]\beta^R_{\kappa_0}$ (provided $j_0>0$) for $\kappa_0=1^{j_1+1}2^{j_1}\dots $, a partition of $R$ into $m+1$ parts and in $J[1]\beta^R_{\kappa_k}$ for each $k=2,\dots ,R$ (provided $j_{k-1}>0$) with $\kappa_k=1^{j_1}\dots (k-1)^{j_{k-1}-1}k^{j_k+1}\dots $, a partition of $R$ into $m$ parts. Using Proposition~\ref{prop:J1Bbold}, we find the coefficient of $\beta^{R-1}_{\kappa}$ in $J[1]\Phi^R$ is
\begin{align*}
\frac{(-1)^{m}}{(m+1)\binom{R}{m+1}}j_0
+\frac{(-1)^{m+1}}{m\binom{R}{m}}\sum_{k=2}^R j_{k-1}=0,
\end{align*}
by using $j_0=R-m$ and $\sum_{k=2}^R j_{k-1}=m$.
\hfill $\Qed$

\begin{proposition}
\label{prop:ZPhiR}
The $1$-point Jacobi function $\Phi^R$ with $R\in \mathbb{N}$ is given by
\begin{align}
Z^J_{V^R_\sigma}\left(\Phi^R;z,\tau\right)
=F_R(z,\tau)\left(\frac{\theta_1(z,\tau)}{\eta(\tau)}\right)^R,
\notag 
\end{align}
\end{proposition}

\noindent \textbf{Proof.}
Using Proposition~\ref{prop:apnpt}, we first note that
\begin{align}
&Z^J_{V^R_\sigma}\left(\beta^{R}_{\pi} ; z,\tau\right)=\binom{R}{j_0j_1\dots j_R}\prod_{\nu=1}^{m}(-1)^{m+1}\widetilde{G}_{k_\nu}(z,\tau) \left(\frac{\theta_1(z,\tau)}{\eta(\tau)}\right)^R.
\label{eq:ZbetaR}
\end{align}
Write $F_R(z,\tau)=Z^J_{V^R_\sigma}\left(\Phi^R;z,\tau\right)\left(\dfrac{\eta(\tau)}{\theta_1(z,\tau)}\right)^R$.
Using \eqref{eq:ZbetaR}, we find that
\begin{align*}
F_R(z,\tau)
&=-\sum_{\pi}\frac{1}{m}\binom{m}{j_1\cdots j_R}
\widetilde{G}_{1}(z,\tau)^{j_1} \cdots \widetilde{G}_{R}(z,\tau)^{j_R}.
\end{align*}
Forming a generating function with parameter $x$ this implies
\begin{align*}
\sum_{R\ge 1} F_R(z,\tau)x^R
&=-\sum_{R\ge 1}\sum_{\pi}\frac{1}{m}\binom{m}{j_1\cdots j_R}
\left(\widetilde{G}_{1}(z,\tau)x\right)^{j_1} \cdots \left(\widetilde{G}_{R}(z,\tau)x^R\right)^{j_R}
\\
&= -\sum_{m\ge 1}\frac{1}{m}\left(\sum_{k\ge 1}\widetilde{G}_{k}(z,\tau)x^k\right)^m
= \text{Log}\left(1-\sum_{k\ge 1}\widetilde{G}_{k}(z,\tau)x^k\right).
\end{align*}
The result follows from \eqref{eq:gen1} and \eqref{eq:P1Gn}.
\hfill $\Qed$

\medskip
We note that $F_R(z,\tau)$ is a meromorphic Jacobi form of weight $R$ for each $R\ge 2$ whereas  $F_1(z,\tau)=P_{1}(z,\tau)$ is quasi-Jacobi.

\medskip
We briefly describe another example. Define  for $R\ge 2$
\begin{align*}
\Psi^R:=\frac{(-1)^R}{R!}\left((R-1)\beta^{R}_{1^R}
+\sum_{k=2}^{R}(-1)^{k+1}(k-1)!\beta^{R}_{1^{R-k}k^1}\right).
\label{eq:PsiR}
\end{align*}

\begin{proposition}
\label{prop:ZPsiR}
For each $R\ge 2$, we have $J[1]\Psi^R=0$ and $\Psi^R$ has $1$-point Jacobi function
\begin{align}
Z^J_{V^R_\sigma}\left(\Psi^R;z,\tau\right)
=K_R(z,\tau) \left(\frac{\theta_1(z,\tau)}{\eta(\tau)}\right)^R,
\notag 
\end{align}
for $K_R(z,\tau)$ given in \eqref{eq:Kn}.
\end{proposition}

\noindent \textbf{Proof.}
Applying Proposition~\ref{prop:J1Bbold} one confirms that $J[1]\Psi^R=0$.
Equation \eqref{eq:ZbetaR} implies that
\begin{align*}
Z^J_{V^R_\sigma}(\Psi^R;z,\tau)\left(\frac{\eta(\tau)}{\theta_1(z,\tau)}\right)^R&=\frac{(R-1)}{R!} \widetilde{G}_{1}(z,\tau)^R
+\sum_{k=2}^{R}\frac{(k-1)!}{R!}\binom{R}{k-1,R-k,1}\widetilde{G}_{1}(z,\tau)^{R-k}\widetilde{G}_{k}(z,\tau)
\\
&=\frac{(R-1)}{R!} \widetilde{G}_{1}(z,\tau)^R
+\sum_{k=2}^{R}\frac{1}{(R-K)!}\widetilde{G}_{1}(z,\tau)^{R-k}\widetilde{G}_{k}(z,\tau)
= K_R(z,\tau).
\end{align*}
\hfill $\Qed$


\end{document}